\documentclass [12pt] {report}
 \usepackage{amsfonts,amsmath,amssymb,amscd,array, euscript}

\newcommand {\Sh} {\mathop{\rm Sh}\nolimits}
\newcommand {\Ch} {\mathop{\rm Ch}\nolimits}

\newcommand {\Cos} {\mathop{\rm Cos}\nolimits}

\newcommand {\Sin} {\mathop{\rm Sin}\nolimits}
\newcommand {\e} {\mathop{\rm e}\nolimits}

\newcommand {\arccosh} {\mathop{\rm arccosh}\nolimits}
\newcommand {\arcsinh} {\mathop{\rm arcsinh}\nolimits}

\newcommand {\D}[2] {\displaystyle\frac{\partial{#1}}{\partial{#2}}}

\newcommand {\al} {\alpha}

\newcommand {\la} {\lambda}

\newcommand {\ka} {\varkappa}
\newcommand {\Ka} {\mbox{\Large$\varkappa$}}

\newcommand {\si} {\sigma}
\newcommand {\Si} {\Sigma}

\newcommand {\ga} {\gamma}

\newcommand {\de} {\delta}
\newcommand {\De} {\Delta}

\newcommand {\prtl} {\partial}
\newcommand {\fr} {\displaystyle\frac}

\newcommand {\be} {\begin{equation}}
\newcommand {\ee} {\end{equation}}
\newcommand {\ba} {\begin{array}}
\newcommand {\ea} {\end{array}}
\newcommand {\bp} {\begin{picture}}
\newcommand {\ep} {\end{picture}}
\newcommand {\bc} {\begin{center}}
\newcommand {\ec} {\end{center}}

\newcommand {\bt} {\begin{tabular}}
\newcommand {\et} {\end{tabular}}
\newcommand {\lf} {\left}
\newcommand {\rg} {\right}

\newcommand {\nin}{\noindent}
\newcommand {\cA} {{\cal A}}

\newcommand {\cC} {{\cal C}}

\newcommand {\cI} {{\cal I}}

\newcommand {\cP} {{\cal P}}

\newcommand {\cR} {{\cal R}}
\newcommand {\cS} {{\cal S}}    \newcommand {\cT} {{\cal T}}
    \newcommand {\cU} {{\cal U}}

\newcommand {\ses} {\medskip}

\newcommand{\bea}{\begin{eqnarray}}
\newcommand{\eea}{\end{eqnarray}}

\def\2#1#2#3{{#1}_{#2}\hspace{0pt}^{#3}}
\def\3#1#2#3#4{{#1}_{#2}\hspace{0pt}^{#3}\hspace{0pt}_{#4}}
\newcounter{sctn}
\def\sec#1.#2\par{\setcounter{sctn}{#1}\setcounter{equation}{0}
                  \noindent{\bf\boldmath#1.#2}\bigskip\par}

\begin {document}

\begin {titlepage}

\vspace{0.1in}

\begin{center}

{\Large \bf   Pseudo-Finsleroid spatial-anisotropic relativistic space}

\end{center}

\vspace{0.3in}

\begin{center}

\vspace{.15in} {\large G.S. Asanov\\} \vspace{.25in}
{\it Division of Theoretical Physics, Moscow State University\\
119992 Moscow, Russia\\
{\rm (}e-mail: asanov@newmail.ru{\rm )}} \vspace{.05in}

\end{center}

\begin{abstract}

\ses

The pseudo-Finsleroid relativistic    metric was constructed upon assuming that
the involved vector field $b_i$ is time-like. In the present paper it is shown
that the metric admits just the alternative counterpart in which the field is space-like.
The entailed pseudo-Finsleroid-spatial framework is systematically described.
We face on  various remarkable properties,
  including the constant  curvature of the associated indicatrix,
  the explicit Hamiltonian function,
  transparent      presentations for the angle    and scalar product.
  The spray coefficients are found to be of a  rather simple       structure.
 The Berwald case is attractively realized.
Interesting  conformal properties  are stemming.

\ses

\end{abstract}

\end{titlepage}

\vskip 1cm

\ses

\ses

\setcounter{sctn}{1} \setcounter{equation}{0}

\nin
  {\bf 1. Introduction: the pseudo-Finsleroid  metric background }

\ses

\ses

The Finslerian  methods propose various interesting tools for addressing
the anisotropy  to the geometry (see [1-7]).
Presence of local spatial anisotropy  in the relativistic space-time
is typical of many pictures appeared in physical applications. Can a preferred spatial
vector field, to be denoted as $b_i(x)$, be {\it geometrically distinguished},
in the sense that the vector $b_i$ can  make traces
on dependence of the fundamental metric function on local directions?
The Finsleroid tools (developed in [8-12])
 can provide us with an adequate method to develop such a subject.

Namely, the simplest standpoint to start from is the stipulation that
 the Finsler space is  specified in accordance with  the condition
that the squared Finslerian metric function
$F^2(x,y)$ possesses  the functional dependence
\be
F^2(x,y) =\Phi \Bigl(g(x), b_i(x), a_{ij}(x),y\Bigr),
\ee
where $g(x)$ is a scalar, $a_{ij}(x)$ is a pseudo-Riemannian metric tensor, and $y$ denote
tangent vectors supported by points $x$ of the underlining manifold.
Next, we require  that
the entailed
Finslerian metric tensor $g_{ij}(x,y)$ be
 of the time-space signature:
\be
{\rm sign}(g_{ij})=(+ - -\dots).
\ee
The pseudo-Riemannian norm $||b||$ of the involved vector $b_i$
should be negative.
Also, when $||b||=-1$, the associated indicatrix should be a space
of constant  curvature.
Finally, sufficient regularity properties are implied.
The respective anisotropic  space
${\mathbf\cA\cR}_{g;c}$ defined in (2.16)
meets hopefully all these requirements.
The space  admits the  norm $||b||$
to be
$||b||=-c$ with $c\in (0,1]$.
The time-like and space-like sectors are separated according to the sign
of the involved  quadratic form $B(x,y)$.
The scalar
\be
h(x)=\sqrt{1+\epsilon \fr{g^2(x)}4},
\ee
where $\epsilon = 1$ in the time-like sector and
$\epsilon = - 1$ in the space-like sector,
plays an important role in the theory developed below.

In Section 2, the     ${\mathbf\cA\cR}$-space as well as  the unimodular
${\mathbf\cU\cA\cR}$-case which uses $c=1$ are rigorously  introduced,
basic notions  are explained, and various key properties are set forth.
The Finslerian metric tensor constructed from
the  function
 $F^2$ given by (2.13)
does inherit from the tensor $a_{ij}(x)$
the time-space signature (2.4), so that (1.2) is valid.
Elucidating the structure of
the respective Cartan tensor and indicatrix curvature tensor
 of the  $ {\mathbf\cA\cR}$-space  results in the remarkable special types
 (2.48) and  (2.50).
       The induced geometry in the tangent spaces of the
                      ${\mathbf\cA\cR}$-space is of the conformally flat type.
When $c=1$,
from (2.43) we can conclude that $X=1/N$, which reduces (2.46) to
the simple representation
 \be
F^2C_hC^h=-\fr{N^2g^2\epsilon}{4},
\ee
which right-hand part is independent of vectors $y$
and, therefore,
 the constancy of the indicatrix curvature is realized (as manifested by
the formulas (2.55)-(2.57)).
At any $c\in (0,1]$, the spray coefficients $G^i$ can  be represented
in the clear explicit form (2.61). The Berwald case thereof is obviously
determined by the conditions $\nabla_j b_i=0$ and $g=const$,
which reduce $G^i$
to the   associated pseudo-Riemannian coefficients
$a^i{}_{nm}(x)y^ny^m$  (where $a^i{}_{nm}(x) $  are the
Christoffel symbols given rise to by the associated pseudo-Riemannian metric
tensor $a_{ij}(x)$).

The   ${\mathbf\cU\cA\cR}$-space angle
 $ \al_{\{x\}}(y_1,y_2)$ formed by two vectors in symmetric way
  is proposed.  The angle
  is   supported  by a point $x\in M$ of the base manifold $M$
(in just the  same sense as in the pseudo-Riemannian geometry)
and is  independent of any vector element
 $y$
 of support.
This gives rise to the notion of the scalar product of the vectors
  $(y_1,y_2)$.
  The  ${\mathbf\cU\cA\cR}$-space coordinates are indicated in terms
of which the metric line element $ds$ can be represented in the transparent forms,
extending the  spherical-coordinate representations handling in the Riemannian geometry.
This $ds$  (see (2.100)) shows that  the  angle
is of the ``canonical'' geometrical meaning:
 $ \al_{\{x\}}(y_1,y_2)$ is  the length of the geodesic curve on the indicatrix,
the curve  joining the points
 at which the entered vectors $(y_1,y_2)$    intersect  the indicatrix.
The success has been predetermined by the observation
that   the squared
{$\mathbf\cA\cR$}-pseudo-Finsleroid metric function $F^2(x,y)$ can conveniently be
represented by the formula
\be
F^2(x,y)=B(x,y) e^{- g(x)\chi(x,y)}
\ee
such that in the $ {\mathbf\cU\cA\cR}$-space
 the scalar  $\chi$ possesses the lucid geometrical meaning
of the azimuthal angle measured from the direction assigned by the input vector $b$.
Indeed, fixing a point $x$ in the background manifold $M$,
a vector $y^i\in T_xM$ is parallel to $b^i(x)$ if and only if
the angle-representation of the $y^i$ corresponds to $\chi=0$.
The property  is also obvious from the
derived formula  (2.76), namely $ \al_{\{x\}}(y,b)=\chi $.
The arisen functions
$\Sh\chi,\,\Ch\chi, \,\Sin\chi,\,\Cos\chi $ can be interpreted as
the required extensions of  the trigonometric functions
to the ${\mathbf\cU\cA\cR}$-space.

    In Section 3,   we confine the attention to the future-time-like sector
        (the respective function $F^2$ introduced by (3.11)-(3.12)
is of the type  proposed  earlier    in the kinematic-study paper [11]).
In Section 4, the space-like sector is attentively evaluated.
The amazing uniformity between structures of the representations of various involved
tensors can be traced, despite the form of the fundamental metric function differs
drastically in these sectors (compare between (3.11)-(3.12) and (4.10)-(4.14)).

In Section 5,   we  fix the tangent space and write down
the respective components of the covariant vector and metric tensor relative to
the orthonormal frame of the input pseudo-Riemannian metric tensor $a_{ij}$,
choosing the $(N-1)$-th component of the tangent vector to be parallel
 to the 1-form $b$. The formulas obtained in this way are convenient in various
 evaluations.

In Section 6,  the possibility of obtaining  the   $\cA\cR$-{\it space  Hamiltonian function} $H$
 in a convenient explicit form is indicated.
When $c=1$,  the function $H^2$  derived is totally similar to the initial  function $F^2$.
In contrast to the later  function $F^2$ which depends primarily on
the contravariant vectors $ y=\{y^i\}$,
the  $H^2$ is the function of the co-vectors $\hat y=\{y_i\}$.
The relationship between the variables $ y=\{y^i\}$ and $\hat y=\{y_i\}$ is assigned
by the simple formula (2.37). The knowledge of the Hamiltonian function makes it possible
to formulate the  Hamilton-Jacobi equation which is of important physical meaning.

In Appendix A,
various calculations with the help of the
 $ {\mathbf\cU\cA\cR}$-space coordinates are systematically presented.

In Appendix B we indicate a straightforward and convenient way to obtain the
${\mathbf\cA\cR}$-space spray coefficients.

  In Appendix C,
 we observe the remarkable conformal property of the ${\mathbf\cU\cA\cR}$-space
that the explicit and simple form can be  proposed for the respective
conformal multiplier,
which proves to be expressed through the fundamental metric function,
according to (C.2).
At the same time, the  $ {\mathbf\cU\cA\cR}$-space
is not conformal to any pseudo-Riemannian space.
Instead,  the  $ {\mathbf\cU\cA\cR}$-space
proves to be metrically isomorphic to the factor-pseudo-Riemannian space ${\cC_g}$
in which the metric tensor reads
\be
t_{mn}(x,\zeta)=p(x,\zeta) a_{mn}(x),
\ee
where  $\zeta\in{\cC_g}$ and the factor is
\be
p(x,\zeta)=\fr1{h^2(x)}|S^2(x,\zeta)|^{(1-h(x))/h(x)},
\ee
with $ S^2(x,\zeta)=a_{ij}(x) \zeta^i\zeta^j$
standing for  the  pseudo-Riemannian metric in the space
$ {\cC_g}$.
The isomorphism is the diffeomorphism, being initiated by the explicit and simple
formulas  (C.1)-(C.4).
There appears the relationship
\be
|F^2(x,y)|^{h(x)}=|S^2(x,\zeta)|,
\ee
which shows that at each point $x$ the Finslerian indicatrix (defined by $|F^2(x,y)|=1$)
is isomorphic to the pseudo-Euclidean sphere $ \cS_x({\cC_g})$
(defined by  $|S^2(x,\zeta)|=1$) in the space $ {\cC_g}$.
The curvature of the pseudo-sphere $ \cS_x({\cC_g})$ in the space  ${\cC_g}$
is not, however, of the unit value, for the metric tensor (1.6)
in the space  ${\cC_g}$ differs from the  pseudo-Riemannian tensor by the factor
$p$.
If we construct on   $ \cS_x({\cC_g})$ the  curvature tensor,
 we readily arrive at the representation (C.29)
which says us that
 in the ${\mathbf\cU\cA\cR}$-space the indicatrix is of the constant curvature
\be
 \cR_{\text{${\mathbf\cU\cA\cR}$-pseudo-Finsleroid  indicatrix} }=- \epsilon h^2.
\ee
This value  agrees with the conclusions (2.55)-(2.57)  obtainable in the traditional
manner when the known tensors of Finsler geometry are systematically evaluated.
Also, the attentive consideration leads to the remarkable equality
\be
\al = \fr1h
 \al_{\text{pseudo-Riemannian}}
 \ee
(see (C.31)),
 where
the right-hand part
 belongs to the sense of the ordinary
pseudo-Riemannian geometry and operates in the
 space $\cC_g$.
The equality (1.10) is a direct implication of (1.6)-(1.8).

To come to the framework with changing  the curvature of the pseudo-Finsleroid indicatrix
along the background manifold of points $x$,
we must allow for  a dependence $g=g(x)$
of the involved pseudo-Finsleroid charge $g$ of $x$, avoiding the stipulation that
the $g$ be a constant.
In this direction, the idea of involutive dependence
(set forth in the previous work [12]),
which permits for  the scalar $g(x)$ to vary just in the direction assigned by the
input preferred vector field $b_i(x)$, seems to be attractive to follow.

 {

\ses

\ses

\setcounter{sctn}{2}            \setcounter{equation}{0}

\nin    {\bf 2.  Basic notions of the  {$\mathbf\cA\cR$}-space}

\ses

\ses

Let $M$ be an $N$-dimensional
$C^{\infty}$
differentiable  manifold, $ T_xM$ denote the tangent space to $M$ at a point $x\in M$,
and $y\in T_xM$  mean tangent vectors.
Suppose we are given on $M$ a pseudo-Riemannian metric ${\cal S}=S(x,y)$.
 Denote by
$\cR_N=(M,{\cal S})$
the obtained $N$-dimensional pseudo-Riemannian space.
Let us also assume that the manifold $M$ admits a non-vanishing space-like 1-form
$ b= b(x,y)$ of the length $c(x)$,            so that
\be
|| b||=-c.
\ee
We shall assume
\be
0<  c(x)\le 1.
\ee

With respect to  natural local coordinates in the space
$\cR_N$
we have the local representations
\be
 b=b_i(x)y^i, \qquad
\cS= a_{ij}(x)y^iy^j,  \qquad   a^{ij}(x)b_i(x)b_j(x)=-c^2(x),
\ee
with a pseudo-Riemannian metric tensor $a_{ij}(x)$, so that
the {\it time-space signature}
\be
{\rm sign}(a_{ij})=(+ - -\dots)
\ee
takes place.
The reciprocity  $a^{in}a_{nj}=\de^i{}_j$ is implied,
 where $\de^i{}_j$ stands for the Kronecker symbol.
The covariant index of the vector $b_i$  will be raised by means of the tensorial rule
$ b^i=a^{ij}b_j,$ which inverse reads $ b_i=a_{ij}b^j.$
We
also  introduce the tensor
\be
r_{ij}(x)~:=a_{ij}(x)+
b_i(x)b_j(x)
\ee
and construct the quadratic form
\be
\ga=r_{ij}y^iy^j,
\ee
in terms of which we introduce the variable
\be
q=\sqrt{|\ga|},
\ee
obtaining
\be
\cS=\ga -    b^2.
\ee

We also introduce on the background manifold $M$
a scalar field $g=g(x)$.
The consideration will be structured by signs of the quadratic form
\be
B(x,y) :=\ga-gbq-b^2\equiv \cS-gbq.
\ee
Namely,  to generalize the pseudo-Riemannian geometry in a desired
pseudo-Finsleroid Finslerian way,
we adapt the consideration to the  following decomposition of the tangent bundle $TM$:
\be
TM=\cT_g^{+}
\cup\Si_g^{+}
\cup\cR_g^{+}
\cup\cR_g^{+0}
\cup\cR\cR_g
\cup\cR_g^{-0}
\cup\cR_g^{-}
\cup\Si_g^{-}
\cup\cT_g^{-},
\ee
which sectors relate to the cases where the tangent vectors $y\in TM$ are, respectively,
time-like, upper-cone isotropic, space-like, lower-cone isotropic, or past-like.
The upperscripts ``$+$'' and  ``$-$'' stand for ``future'' and ``past'', respectively.
We take the sums
\be
\cT_g^{+}=\cT_g^{+r}\cup\cT_g^{+l}, \quad   \Si_g^{+}=\Si_g^{+r}\cup\Si_g^{+l},
\qquad
\Si_g^{-}=\Si_g^{-r}\cup\Si_g^{-l}, \quad   \cT_g^{-}=\cT_g^{-r}\cup\cT_g^{-l},
\ee
\ses
\be
\cR_g^{+}=\cR_g^{+r}\cup\cR_g^{+l}, \qquad
\cR\cR_g=\cR\cR_g^{r}\cup\cR\cR_g^{l},   \qquad
\cR_g^{-}=\cR_g^{-r}\cup\cR_g^{-l},
\ee
where
$r$ means {\it in the direction of the vector} $b^i$
\, (= to the right), and $l$ means {\it opposite to the direction of} $b^i$
\, (= to the left).

We propose

\ses\ses

{\large Definition}. The squared
{$\mathbf\cA\cR$}-{\it pseudo-Finsleroid metric function} $F^2(x,y)$ is given
by the formula
\be
F^2(x,y)~:=B(x,y)J^2(x,y)
\ee
with the function $J$ taken from (3.12) and (4.11).

\ses

\ses

The positive  (not absolute) homogeneity  holds:
 $F^2(x,\la y)=\la^2 F^2(x,y)$ for any $\la >0$   and all admissible $(x,y)$.

We introduce the {\it indicator} $\epsilon$:
\be
\epsilon=1, \,\, \text{if vector $y$ is time-like}, \quad   \text{and} \quad
\epsilon=-1, \,\, \text{if vector $y$ is space-like}.
\ee

 We can  obtain from the function $F^2$ the  distinguished Finslerian tensors,
 and first of all
the covariant tangent vector $\hat y=\{y_i\}$,
the  Finslerian metric tensor $\{g_{ij}\}$
together with the contravariant tensor $\{g^{ij}\}$ defined by the reciprocity conditions
$g_{ij}g^{jk}=\de^k_i$, and the  angular metric  tensor
$\{h_{ij}\}$, by making  use of the  conventional  Finslerian  rules
 in succession:
 \be
y_i :=\fr12\D{F^2}{y^i}, \qquad
g_{ij} :
=
\fr12\,
\fr{\prtl^2F^2}{\prtl y^i\prtl y^j}
=\fr{\prtl y_i}{\prtl y^j}, \qquad
h_{ij} := g_{ij} -\fr1{F^2}y_iy_j.
\ee

Under these conditions, we introduce

\ses

\ses

 {\large  Definition}.  The arisen  space
\be
{\mathbf\cA\cR}_{g;c} :=\{\cR_{N};\, TM;\, y\in TM;\,b_i(x);\,g(x);\,F^2(x,y);\,g_{ij}(x,y)\}
\ee
is called the ${\mathbf\cA\cR}$-{\it pseudo-Finsleroid space}.

\ses

\ses

 {\large  Definition}. The space $\cR_N=(M,\,{\cal S})$ entering the above definition
 is called the {\it associated  pseudo-Riemannian space}.

\ses

 \ses

 {\large  Definition}.  Within  any tangent space $T_xM$, the
  function $F^2(x,y)$
 produces
the ${\mathbf\cA\cR}$-{\it  pseudo-Finsleroid indicatrix}:
 \be
 \cI{\mathbf\cA\cR}_{g;c\, \{x\}} :=\{y\in \cI{\mathbf\cA\cR}_{g;c\, \{x\}} :
  y\in T_xM, F^2(x,y)=1\},      \quad \text{time-like} \,\, \{y\};
  \ee
\ses
 \be
 \cI{\mathbf\cA\cR}_{g;c\, \{x\}} :=\{y\in \cI{\mathbf\cA\cR}_{g;c\, \{x\}} :
  y\in T_xM, F^2(x,y)=0\},     \quad \text{isotropic} \, \, \{y\};
  \ee
\ses
 \be
 \cI{\mathbf\cA\cR}_{g;c\, \{x\}} :=\{y\in \cI{\mathbf\cA\cR}_{g;c\, \{x\}} :
  y\in T_xM, F^2(x,y)=-1\}, \quad \text{space-like} \,\, \{y\}.
  \ee

\ses

 {\large  Definition}.
 $ \cI{\mathbf\cA\cR}_{g;c\, \{x\}} \subset T_xM$ is the boundary of
  the ${\mathbf\cA\cR}$-{\it  pseudo-Finsleroid}
    $ \cP{\mathbf\cA\cR}_{g;c\, \{x\}} \subset T_xM$.

 \ses

 {\large  Definition}. The scalar $g(x)$ is called
the ${\mathbf\cA\cR}$-{\it pseudo-Finsleroid charge}.
The 1-form $b=b_i(x)y^i$ is called the  ${\mathbf\cA\cR}$-{\it pseudo-Finsleroid-axis}
 1-{\it form}.

{

In the time-like sectors the quadratic form (2.9) reads merely
$$
B =q^2-gbq-b^2>0,
$$
and, therefore,
is of the positive discriminant
$$
D_{\{B\}}=4h^2>0, \qquad   \epsilon=1,
$$
with  $h=\sqrt{1+(g^2/4)}$.
Alternatively,
in the space-like sectors, the formula (2.7) yields $\ga=-q^2$
and  the $B$ of  (2.9) takes on the form
$$
B =-(q^2+gbq+b^2) < 0,
$$
which   discriminant  is negative:
$$
D_{\{B\}}=-4h^2<0, \qquad   \epsilon=-1,
$$
where this time $h=\sqrt{1-(g^2/4)}$.
This distinction is the reason why the form of the pseudo-Finsleroid  function
$F^2$
is essentially different  (see (3.11)-(3.12) and (4.10)-(4.13))
in the time-like and space-like sectors.

{

In the limit $g\to 0$,
the definition (2.9) degenerates to the
 quadratic form  of the input pseudo-Riemannian metric tensor
 $a_{ij}(x)$:
$$
B|_{_{g=0}}=
a_{ij}(x)y^iy^j
\equiv \cS(x,y).
$$
We have also
\be
J(x,y)|_{_{g=0}}=1,
\qquad
F^2(x,y)|_{_{g=0}}=
 \cS(x,y),
\ee
and
\be
a_{ij}(x)=g_{ij}(x,y)\bigl|_{g=0}\bigr.    .
\ee

\ses

\ses

 {\large  Definition}. The space ${\mathbf\cA\cR_{g;c}}$ is  {\it unimodular},
if $c=1$:
\be
{\mathbf\cU\cA\cR_g} :={\mathbf\cA\cR_{g;c=1}}.
\ee

\ses

\ses

The derivatives
\be
\D Jb=-\fr g2\fr qBJ, \qquad
\D Jq=\fr g2\fr bBJ
\ee
can readily be obtained from (3.12) and (4.11).

In many cases it is convenient to use the variables
\be
u_i~:=a_{ij}y^j,
\qquad
v^i~:=y^i+bb^i, \qquad v_m~:=u_m+bb_m=r_{mn}y^n     \equiv a_{mn}v^n.
\ee
The identities
\be
u_iv^i=v_iy^i=\epsilon q^2,
 \qquad
v_ib^i=v^ib_i=(1-c^2)b,   \qquad
 r_{ij}b^j=(1-c^2)b_i,
\ee
\ses
\be
v_iv^i=\epsilon q^2+(1-c^2)b^2,  \qquad
\D b{y^i}=b_i, \qquad \D q{y^i}=\epsilon   \fr{v_i}q,
\qquad \D {(q/b)}{y^i}=\fr q{b^2}
e_i
\ee
can readily be verified,
where
\be
e_i=-b_i+ \epsilon   \fr b{q^2}v_i
\ee
is the vector showing the property
$
e_iy^i=0.
$
 Particularly, the vector enters the equality
\be
\D{\lf(J^2\rg)} {y^k}=  \fr{gq}B J^2      e_k.
\ee

{

By performing required direct calculation,
we find
the representations
\be
y_i=\bigl(v_i-(b+gq)b_i\bigr) J^2
\ee
and
\be
g_{ij}=
\Biggl[a_{ij}
-\fr g{B}\lf ( -   q(b+gq)b_ib_j+q(b_iv_j+b_jv_i)-
\epsilon
b\fr{v_iv_j}q\rg)\Biggr]J^2,
\ee
together with the reciprocal (contravariant) components
\be
g^{ij}=
\Biggl[a^{ij}
+\fr g{B}\lf(-bqb^ib^j   +    q(b^iv^j+b^jv^i)  -
\epsilon
 (b+gc^2q)\fr{v^iv^j}{\nu}\rg)
\Biggr]\fr1{J^2}.
\ee
The determinant of the metric tensor  (2.30)
can straightforwardly be evaluated, yielding
\be
\det(g_{ij})= \fr{\nu}q   J^{2N} \det(a_{ij}),
\ee
where
\be
\nu=q-  \epsilon   (1-c^2)gb.
\ee

\ses

It is useful to note that
\be
b(b+gc^2q)=-B+   \epsilon   q\nu,
\ee
\ses
\be
g_{ij}b^j=
\Biggl[b_i
-\fr g{B}\lf (
   q(b+gc^2q)b_i
                 -\fr 1q\bigl(c^2q^2  + \epsilon (1-c^2)b^2\bigr)v_i\rg)
\Biggr]J^2,
\ee
\ses
and
\be
g_{ij}v^j=      \fr{\nu}q
\Biggl[
   (\epsilon q^2-b^2)v_i
-\epsilon gq^3b_i
\Biggr]\fr1BJ^2.
\ee

{

In terms of the set $\{b_i,u_i=a_{ij}y^j \}$, we obtain the alternative representations
\be
y_i=\bigl(u_i - gqb_i\bigr)J^2
\ee
\ses
and
\be
g_{ij}=
\biggl[a_{ij}
+\fr g{B}\Bigl((gq^2-\fr{b(q^2-\epsilon b^2)}q)b_ib_j  + \epsilon  \fr bqu_iu_j
-
\fr{ q^2 -\epsilon b^2}q(b_iu_j+b_ju_i)\Bigr)\biggr]J^2.
\ee
We get also
\be
g^{ij}=
\biggl[a^{ij}+ \epsilon  \fr g{\nu}(bb^ib^j+b^iy^j+b^jy^i)
- \epsilon \fr g{B\nu}(b+gc^2q)y^iy^j
\biggr]\fr 1{J^2}.
\ee
With the help of (2.37) we can transform (2.38) to
\be
g_{ij}=
\biggl[a_{ij}
+  \epsilon    \fr g{q}
\Bigl(
-(b+gq)b_ib_j
+   \fr1B b   \fr1{J^2} \fr1{J^2}  y_iy_j
-  \fr 1{J^2}(b_iy_j+b_jy_i)\Bigr)
\biggr]J^2.
\ee

To raise the index, it is convenient to apply the rules
$$
g^{ij}b_j=\fr1{F^2}\Biggl[(B+gbq)b^i
-\fr {\epsilon g}{\nu}\lf(c^2B+b(b+gc^2q)\rg)v^i\Biggr]
$$
\ses
and,
for any co-vector $t_j$,  from (2.31) we obtain
$$
g^{ij}t_j=
\Biggl[Ba^{ij}t_j
+gq(yt)b^i+\fr {\epsilon g}{\nu}
\Bigl(B(bt)-(b+gc^2q)(yt)\Bigr)v^i
\Biggr]\fr 1{F^2},
$$
where
$(yt)=y^jt_j$  and  $(bt)=b^jt_j$.

{

From the determinant  value (2.32) we can explicate the vector
\be
C_i=
   \D{\ln\lf(\sqrt{|\det(g_{mn})|}\rg)}{y^i},
  \ee
 obtaining
\be
C_i= g\fr {1}  {2B}
\fr qX
e_i,
\ee
where
\be
\fr1X=    N+  \epsilon  \fr{(1-c^2)B}{q\nu}.
\ee
Also,
\be
C^i=g\fr 1{2F^2}
\fr qX
\lf(-b^i  + \epsilon  \fr 1{q\nu}(b+gc^2q)v^i\rg),
\ee
or
alternatively,
\be
C^i= \epsilon    g\fr 1{2F^2}
\fr 1{X\nu}
\Bigl(-Bb^i  +(b+gc^2q)y^i\Bigr).
\ee
We can evaluate  the contraction
\be
C^hC_h=
  - \epsilon  \fr{g^2}4
\fr1{F^2X^2}\lf(N+1-\fr1X\rg).
\ee

The    Cartan tensor
\be
C_{ijk}~ := \fr 12\D{g_{ij}}{y^k},
 \ee
when evaluated with the help of the components $g_{ij}$ given by (2.30),
is representable in the form
\be
C_{ijk}= X
 \Biggl[
C_ih_{jk}  +C_jh_{ik}  +C_kh_{ij}
-\lf(N+1-\fr1X\rg)
\fr1{C_hC^h}C_iC_jC_k
\Biggr].
\ee

\ses

{

The curvature of  indicatrix
is described by
the tensor
\be
\hat R_i{}^j{}_{mn} := \3Chjm\3Cihn-\3Chjn\3Cihm.
\ee
Inserting  here (2.48)
results in
$$
\hat R_{ijmn}=-(C_hC^h)X^2\bigl(h_{im}h_{jn}-h_{in}h_{jm}\bigr)
$$

\ses

\be
+X^2\lf(N-\fr1X\rg)
\Bigl(C_iC_mh_{jn}-C_iC_nh_{jm}+C_jC_nh_{im}-C_jC_mh_{in}\Bigr).
\ee
Contracting this tensor yields
\be
\hat R_{ijmn}g^{jn}=-(C_hC^h)X^2(N-2)h_{im}
+X^2\lf(N-\fr1X\rg)
\Bigl((N-3)C_iC_m+C_jC^jh_{im}\Bigr)
\ee
\ses
and
\be
\hat R_{ijmn}g^{jn}g^{im}=
-(C_hC^h)X^2(N-2)(N-1)
+X^2\lf(N-\fr1X\rg)
(C_hC^h)(2N-4).
\ee
From these formulas it can readily be concluded that the
tangent pseudo-Riemannian space
 is conformally flat. For instance, taking the dimensions
$N>3$, we can evoke the conformal Weyl tensor
$W_{ijmn}$ and use the representations (2.50)-(2.52). We straightforwardly obtain
$W_{ijmn}=0.$

{

In  the unimodular space ${\mathbf\cU\cA\cR}_g$
defined by (2.22), the equality  $c=1$
entails $X=1/N$ (see (2.43)), so that
that the representation (2.48) reduces to
\be
C_{ijk}=\fr1N\lf(h_{ij}C_k+h_{ik}C_j+h_{jk}C_i-\fr1{C_hC^h}C_iC_jC_k\rg)
\ee
with (1.2) being valid.
The curvature tensor representation (2.50) reduces to merely
\be
F^2\hat R_{ijmn}=\fr{\epsilon g^2}4  (h_{im}h_{jn}-h_{in}h_{jm}),
\ee
which entails the following formulas to characterize the value $\cR$ of the indicatrix
curvature:
\be
 \cR_{\text{${\mathbf\cU\cA\cR}$-pseudo-Finsleroid  indicatrix} }=-\lf(1+\fr14g^2\rg) \le -1,
 \quad when \quad \epsilon = 1,
\ee
 and
\be
 \cR_{\text{${\mathbf\cU\cA\cR}$-pseudo-Finsleroid  indicatrix} }=\lf(1-\fr14g^2\rg) \le 1,
 \quad when \quad \epsilon = -1.
\ee
We have
\be
\cR_{\text{${\mathbf\cU\cA\cR}$-pseudo-Finsleroid indicatrix} }\stackrel{g\to 0}{\Longrightarrow}
\cR_{\text{pseudo-Euclidean sphere}}.
\ee
{The ${\mathbf\cU\cA\cR}$-{\it  pseudo-Finsleroid indicatrix  of the time-like case is a space
 of constant  negative curvature} (according to (2.55)).
 The positive curvature value (2.56) is obtained in the space-like case.

{

With the expression $g^{ij}t_j$ indicated before (2.41), it follows that
$$
F^2(b+gc^2q)g^{ij}t_j=
\Bigl[Ba^{ij}t_j
+gq(yt)b^i
\Bigr](b+gc^2q)
$$

\ses

$$
+\fr {\epsilon g}{\nu}
\Bigl(B(bt)-(b+gc^2q)(yt)\Bigr)(b+gc^2q)v^i.
$$
Since
$$
\fr{2F^2X}gC^i +q b^i  =  \fr{\epsilon}{\nu}   (b+gc^2q)v^i
$$
(see (2.44)), we can write
$$
F^2(b+gc^2q)g^{ij}t_j=
Ba^{ij}t_j
(b+gc^2q)
+
(bt)  q
gBb^i
$$

\ses

$$
+
2\Bigl(B(bt)-(b+gc^2q)(yt)\Bigr) F^2XC^i.
$$
Use here
$$
gBb^i =g(b+gc^2q)y^i
-\fr{2F^2X\nu}{\epsilon }
C^i
$$
(see (2.45))
and the equality
$B=-b(b+gc^2q) +   \epsilon   q\nu
$
(see (2.34))   to   get
$$
F^2g^{ij}t_j=
Ba^{ij}t_j
+
(bt) g q
y^i
-2bF^2X(bt)C^i
-
2(yt) F^2XC^i.
$$
\ses
In this way we arrive at the expansion
$$
F^2g^{ij}t_j=
B\Bigl(a^{ij}t_j
+\fr1{c^2}(bt)b^i\Bigr)
+2 \fr q{gc^2}  F^2X(bt)C^i
$$

\ses

\be
-
2\Bigl(\fr1{c^2}b(bt)+(yt)\Bigr) F^2XC^i
-\fr1{c^2}b(bt)y^i.
\ee

{

We use
 the pseudo-Riemannian covariant derivative
\be
\nabla_ib_j~:=\partial_ib_j-b_ka^k{}_{ij},
\ee
where
\be
a^k{}_{ij}~:=\fr12a^{kn}(\prtl_ja_{ni}+\prtl_ia_{nj}-\prtl_na_{ji})
\ee
are the
Christoffel symbols given rise to by the associated pseudo-Riemannian metric ${\cal S}$.

Attentive direct calculations
(see Appendix B)
 of    the induced spray coefficients
$ G^i =\ga^i{}_{nm}y^ny^m$,
where
$\ga^i{}_{nm}$ denote the associated Finslerian Christoffel symbols,
can be used to arrive at the following result.
\ses
\ses

PROPOSITION.  {\it In the
${\cA\cR}$-space
the  spray coefficients $ G^i$
can  explicitly be written in the form
\be
G^i=
-\epsilon \fr  g{\nu}
\Bigl(
y^jy^h\nabla_jb_h
-gqb^jf_j\Bigr)
v^i
+gqf^i
 +E^i
 +a^i{}_{nm}y^ny^m,
 \ee
 with the  notation
$v^i=y^i+bb^i  $
and
\be
f^i=f^i{}_ny^n,\quad
f^i{}_n=a^{ik}f_{kn}, \quad
f_{mn}=
\nabla_mb_n-\nabla_nb_m
\equiv \D{ b_n}{x^m}-\D {b_m}{x^n}.
\ee
The coefficients $E^i$ involve the gradients
 $g_h=\partial g/\partial x^h$
and  can be taken as
\be
E^i = g^{ih}\fr{\partial y_{h}}{\partial g} (yg)
-\fr12 {\bar M}F^2g_hg^{ih},
\ee
where $X$ is the function given in}
 (2.43), $(yg)=y^hg_h$,
 {\it  and
the function ${\bar M}$  is defined by the equality
\be
\D{F^2} g= {\bar M}F^2.
\ee
}

\ses

\ses

Applying the  formula (2.58) to the case $t_i=g_i$ yields
$$
F^2 g_hg^{ih} =
B\Bigl(g^i
+\fr1{c^2}(bg)b^i\Bigr)
+2 \fr q{gc^2}  F^2X(bg)C^i
$$

\ses

\be
-
2\Bigl(\fr1{c^2}b(bg)+(yg)\Bigr) F^2XC^i
-\fr1{c^2}b(bg)y^i,
\ee
where $(bg)=b^hg_h$.

{

In (2.61),
the difference
$
G^i-a^i{}_{nm}y^ny^m
$
involves  the  crucial terms linear in the covariant derivative
 $ \nabla_jb_h.
 $
When $g=const$,  we have $E^i=0$ and   the  metric function $F$ does not
enter the right-hand side of (2.61),
 the only trace of the function $F$ in the
spray coefficients (2.61) being the occurrence of the pseudo-Finsleroid charge
 $g$ in the right-hand side.


Now we perform differentiation
with respect to the pseudo-Finsleroid parameter $g$, obtaining
 \be
  \D hg= \epsilon
  \fr14 G, \qquad  \D Gg= \fr1{h^3}.
  \ee
\ses
From (4.12) and (4.13) we find
\be
\D fg= -\fr1{2h} +  \fr{b}{-B}  \Bigl(\fr14G  q +   \fr1{2h}b\Bigr),
\quad \epsilon =-1,
\ee
and
 $$
  \D {F^2}g=  -  bqJ^2 - \fr1{h^3}fF^2+
G \Biggl[\fr1{2h}-\fr{b}{-B}  \Bigl(\fr14G  q +   \fr1{2h}b\Bigr)\Biggr]  F^2
={\bar M}  F^2,
\quad \epsilon =-1,
$$
where
$$
{\bar M}=
- \fr{b q}B - \fr1{h^3}f+
\fr12\fr{G}{-Bh} ( q^2+\frac12 gb q),
\quad \epsilon =-1,
$$
\ses
or
\be
{\bar M}=
 - \fr1{h^3}f+
\fr12\fr{G}{-Bh}   q^2+  \frac1{-Bh^2} b q,
\quad \epsilon =-1.
\ee
We find the vector $\bar M_i=\partial \bar M/\partial y^i$:
\be
{\bar M}_i=\fr{4  q^2}{gB} XC_i,
\quad \epsilon =-1.
\ee
\ses

{

\nin
Differentiating (2.64) with respect to $y^i$  and using (2.69) just yield
\be
 \D {{y_i}}g=
\bar My_i
+\fr{2  q^2}{g}J^2 XC_i,
\quad \epsilon =-1.
\ee
Inserting this result  in  (2.63) yields
$$
E^i =
{\bar M} (yg)y^i  +  F^2 \fr{2 q^2 }{gB}(yg) X C^i
-\fr12 {\bar M}F^2g_hg^{ih},
\quad \epsilon =-1.
$$

Similar expressions
 can be derived in the case
$ \epsilon =1$  (see (3.12)-(3.13)).

{

The coefficients $G^i$ determine the geodesic equation
$$
\fr{d^2x^i}{ds^2}+G^i\lf(x,\fr{dx}{ds}\rg)=0,
$$
where
$ds=\sqrt{|F^2(x,dx)|}$.

\ses

Also, it is  possible to   transfer to the present theory the concept
of angle that was introduced in the previous work [8]
 (dealt with the case when the involved vector field $b_i(x)$ is time-like).
Indeed, whenever $c=1$ (whence $||b||=-1$),
  any two nonzero tangent vectors
  $y_1,y_2\in T_xM$ of a fixed tangent space
form the ${\mathbf\cU\cA\cR}$-{\it angle}
 \be
 \al_{\{x\}}(y_1,y_2)~: =
 \fr1h\arccosh
 \fr{ h^2\langle y_1,y_2\rangle _{\{x\}}^{\{r\}}  -  A(x,y_1)A(x,y_2)  }
 { \sqrt{|B(x,y_1)|}\,\sqrt{|B(x,y_2)|} },
 \ee
in the future-time-like sector,
and
  \be
 \al_{\{x\}}(y_1,y_2)~: =
 \fr1h\arccos
 \fr{ A(x,y_1)A(x,y_2)- h^2\langle y_1,y_2\rangle _{\{x\}}^{\{r\}}    }
 { \sqrt{|B(x,y_1)|}\,\sqrt{|B(x,y_2)|} },
 \ee
in the space-like sector,
where $\langle y_1,y_2\rangle _{\{x\}}^{\{r\}}=r_{ij}(x)y_1^iy_2^j$
and
$A=b+(1/2)gq$.

At equal vectors, the zero-value
\be
 \al_{\{x\}}(y_1,y_2)=0
\ee
 takes place.
The zero-degree homogeneity
\be
 \al_{\{x\}}(y_1,ky_2)=
 \al_{\{x\}}(ky_1,y_2)= \al_{\{x\}}(y_1,y_2),
 \qquad k>0,\, \forall y_1,y_2,
\ee
holds.
The proposal   obviously exhibits the symmetry:
\be
 \al_{\{x\}}(y_1,y_2)= \al_{\{x\}}(y_2,y_1).
\ee
The  angle $\al_{\{x\}}(y_1,y_2)$
is   supported  by a point $x\in M$ of the base manifold $M$
(in just the  same sense as in the pseudo-Riemannian geometry),
and is  independent of any vector element
 $y$
 of support.

If, fixing a point $x$,  we consider the angle
$ \al_{\{x\}}(y,b)$ formed by a vector $y\in T_xM$ with
 the input characteristic vector $b^i(x)$,
taking into account that the nullification $q=0$ and the equality
$A=-1$ take place at  $y=b$
 whenever the
  ${\mathbf\cU\cA\cR}$-pseudo-Finsleroid space is used,
from (2.72) we get in the space-like sector
the respective value to be
 \be
 \al_{\{x\}}(y,b)~: = \fr1h\arccos \fr{ A(x,y)}
 {\sqrt{|B(x,y)|}}.
 \ee

{

In the Euclidean and Riemannian geometries,
an important role is played by the spherical coordinates.
Their use enables to conveniently represent vectors,
evaluate squares and volumes,
study curvature of surfaces,
in many cases simplify consideration
and solve rigorously equations,
and also introduce and use
various trigonometric functions.
Similar coordinates are available
in the pseudo-Euclidean geometry and applications founded  upon such the geometry.
In the context of the present theory, such coordinates can meaningfully be extended.
The key observation thereto is the appropriate choice of the angle $\chi$
 as follows:
\be
\chi =\fr1hf
\ee
in agreement with
\be
J=e^{-\frac12g\chi}
\ee
(see (4.11) and formulas below  (3.12)).

Accordingly, by fixing the tangent space according to Section 5,
and choosing  the four-dimensional case $N=4$,
the
${\mathbf\cU\cA\cR}$-{\it space coordinates}
$\{z^p\}$
are given by
\be
z^0=\sqrt{|F^2|}, \qquad  z^1=\eta, \qquad z^2=\phi, \qquad z^3=\chi,
\ee
where  $\phi$ is the {\it polar angle} in the $R^1\times R^2$-plane,
$\chi$ plays the role of the   pseudo-Finsleroid azimuthal angle measured from the direction
of the input vector $b^i$, and $\eta$ serves to measure the  time-components.
The indices $p,q,...$ will be specified over the range 0,1,2,3.

\ses\ses

For the vectors $\{R^p\}$  (see (5.4)-(5.7))
it proves possible  to construct
 the representation
\be
R^p=R^p(g;z^q)
\ee
which
  possesses the {\it invariance property}
\be
F^2(g;R^p(g;z^q))= (z^0)^2\epsilon.
\ee
We shall transform the Finslerian metric
tensor $g_{pq}(g;R)$ given by the components
(5.17)-(5.20) to such coordinates, obtaining the tensor
\be
A_{rs}(g;z) ~:=
g_{pq}(g;R)\D{R^p}{z^r}\D{R^q}{z^s}
\ee
which is of the diagonal structure.

{

In   the {\it future-time-like sector}, so that $\epsilon=1$,
we have $F>0$ and can take $z^0=F$.
The  following representations are appropriate:
\ses
\be
R^0=F\cosh\eta \Ch\chi,
 \ee
\ses
\be
R^1=F\sinh\eta  \Ch\chi\cos\phi,
\qquad
R^2=F\sinh\eta \Ch\chi \sin\phi,
\qquad
R^3=F \Sh\chi ,
\ee
with
\be
\Ch\chi=\fr1{Jh} \cosh f, \qquad
\Sh\chi=
\fr1J\lf(\sinh f-\fr G2\cosh f\rg),
\ee
entailing
\be
  q=F\Ch\chi, \qquad
   b+\fr12gq =\fr FJ\sinh f.
\ee

In the {\it space-like sector}, when $  \epsilon=-1$ is used,
we introduce the function
\be
K=\sqrt{|F^2|}>0, \quad \text{so that} \quad z^0=K,
\ee
and set forth
the representations
\be
R^0=K\sinh\eta \Sin\chi,
 \ee
\ses
\be
R^1=K\cosh\eta  \Sin\chi\cos\phi,
\qquad
R^2=K\cosh\eta \Sin\chi \sin\phi,
\qquad
R^3=K \Cos\chi,
\ee
with
\be
\Sin\chi=\fr1{Jh} \sin f, \qquad
\Cos\chi=
\fr1J\lf(\cos f-\fr G2\sin f\rg),
\ee
which entails
\be
q=K\Sin\chi,  \qquad
 b+\fr12gq =\fr KJ\cos f.
\ee

The arisen functions
$\Sh\chi,\,\Ch\chi, \,\Sin\chi,\,\Cos\chi $ can be interpreted as
the required extensions of  the trigonometric functions
to the ${\mathbf\cU\cA\cR}$-space.

When these representations are applied to  the ${\mathbf\cU\cA\cR}$-angle (2.71)-(2.72),
 the following result is obtained:
\be
 \al_{\{x\}}(y_1,y_2)    = \fr1h\arccosh \tau_{12}
 \ee
with
\be
\tau_{12}=
\cosh(f_2-f_1)
+\Omega_{12}\,\cosh f_1\cosh f_2
\ee
\ses
and
\be
\Omega_{12}=\cosh(\eta_2-\eta_1) -1+\bigl(1-\cos(\phi_2-\phi_1)\bigr)\,\sinh \eta_1 \sinh\eta_2
\ee
 in the         future-time-like sector;
alternatively,
\be
 \al_{\{x\}}(y_1,y_2)    = \fr1h\arccos \tau_{12}
 \ee
with
\be
\tau_{12}=
\cos(f_2-f_1)+\Omega_{12}\,\sin f_1\sin f_2
\ee
\ses
and
\be
\Omega_{12}=\cosh(\eta_2-\eta_1) -1-\bigl(1-\cos(\phi_2-\phi_1)\bigr)\,\cosh \eta_1 \cosh\eta_2
\ee
in the space-like  sector.

{

From (2.82) we get
the diagonal tensor, so that $A_{01}=A_{02}=A_{03}=A_{12}= A_{13}= A_{23}=0$,
and the squared linear element  $ds^2=\epsilon A_{rs}dz^rdz^s$
which reads
\be
(ds)^2=(dz^0)^2-(z^0)^2\Bigl[
(d\chi)^2+\fr1{h^2}\cosh^2 f
\bigl(\sinh^2\eta\,(d\phi)^2+(d\eta)^2\bigr)
\Bigr]
\ee
\ses
in the         future-time-like sector,
and
\be
(ds)^2=(dz^0)^2 +
(z^0)^2\Bigl[(d\chi)^2-
\fr1{h^2}\sin^2 f
\bigl((d\eta)^2-\cosh^2\eta\,(d\phi)^2\bigr)
\Bigr]
\ee
in the space-like  sector.

If we consider the infinitesimal version of the representations (2.92)-(2.97),
putting
$\eta_2-\eta_1=d\eta, \, \chi_2-\chi_1=d\chi,$ and $ \phi_2-\phi_1=d\phi $,
then from (2.98) and (2.99) we can directly conclude that
\be
(ds)^2=(dz^0)^2-(z^0)^2(d\al)^2,
\ee
where $d\al$ is the value issued from (2.92) and (2.95).
The formula (2.100)
is remarkable because tells us that,
with the angle $\al$ defined by
(2.71)-(2.72),
 $d\al$ {\it is the arc-length on the indicatrix}, for along the indicatrix
 we have $z^0=1$.

With the angle
(2.71)-(2.72),
we can naturally propose
  the ${\mathbf\cU\cA\cR}$-{\it scalar product}
 \be
\langle y_1,y_2\rangle _{\{x\}}~:=F(x,y_1)F(x,y_2) \cosh\Bigl( \al_{\{x\}}(y_1,y_2)\Bigr),
 \ee
in the future-time-like sector,  and
 \be
\langle y_1,y_2\rangle _{\{x\}}~:=K(x,y_1)K(x,y_2) \cos\Bigl( \al_{\{x\}}(y_1,y_2)\Bigr),
 \ee
in the space-like sector,  where $K=\sqrt{|F^2|}>0$.
At equal vectors, the reduction
\be
\langle y,y\rangle_{\{x\}}=F^2(x,y)
\ee
 takes place, that is, the two-vector scalar product proposed
reduces  exactly  to  the squared  ${\mathbf\cU\cA\cR}$-Finsler metric function.
The homogeneity
\be
\langle y_1,ky_2\rangle_x=k\langle y_1,y_2\rangle_x,
\qquad
\langle k y_1,y_2\rangle_x=k\langle y_1,y_2\rangle_x,
\qquad k>0,\, \forall y_1,y_2,
\ee
holds.
The proposal   obviously exhibits the symmetry:
\be
\langle y_1,y_2 \rangle_{\{x\}}=  \langle y_2,y_1 \rangle _{\{x\}}.
\ee

{

\ses

\ses

\setcounter{sctn}{3}
\setcounter{equation}{0}

\nin
{\bf 3.  Future-time-like sector of the   space ${\cA\cR}_{g;c}$}

\ses

\ses

Assuming
\be
-\infty<g(x)<\infty,
\ee
we construct the scalar
\be
h(x)=\sqrt{1+\fr14(g(x))^2}.
\ee
In the {\it future-time-like sector}
\be
y\in\cT_g^{+},  \qquad            \cT_g^{+}=\cT_g^{+r}\cup\cT_g^{+l}
\ee
with
\be
\cT_g^{+r}=\Bigl(y\in \cT_g^{+r}:~y\in T_xM, \,  b\ge0,\,  q>-g_-b\Bigr)
\ee
and
\be
\cT_g^{+l}=\Bigl(y\in \cT_g^{+l}:~y\in T_xM, \,  b\le0,\, q> -  g_+b\Bigr),
\ee
where
$$
g_+=-\fr12g+h,  \quad  g_-=-\fr12g-h,
$$
 we have
\be
r_{ij}y^iy^j>0
\ee
and can write the variable (2.7) merely as
\be
q=\sqrt{r_{ij}(x)y^iy^j}.
\ee
The quadratic form (2.9) reads now
\be
B =q^2-gbq-b^2   =   (b-g_-q)  (g_+q-b)>0.
\ee
Notice that
\be
b-g_-q>0,   \qquad g_+q-b>0
\ee
throughout the sector (3.3).
The form $B$ can also conveniently be represented as
\be
B=h^2q^2-\lf(b+\fr12gq\rg)^2.
\ee

{

Throughout the sector (3.3),
we take the metric function (2.13) as follows:
\be
F(x,y)~:=\sqrt{B(x,y)} \,  J(x,y)
\equiv
(b-g_-q)^{G_+/2}(g_+q-b)^{-G_-/2}
\ee
with
\be
J(x,y)=
\lf(
\fr{b-g_-q}{g_+q-b}
\rg)^{-G/4},
\ee
where
$G=g/h$, \,
$G_+=g_+/h$, and $G_-=g_-/h$.

\ses

\ses

Introducing the function
$$
f=\fr12\ln\lf( \fr{b-g_-q}{g_+q-b}  \rg),
$$
we obtain the representation
$$
J(x,y)=\e^{-\frac12G(x)f(x,y)}
$$
of the type (4.11).
We can write
$$
\cosh f=\fr{hq}{\sqrt B}, \qquad \sinh f=\fr{b+\fr12gq}{\sqrt B}.
$$

With the function
$$ L=q-\fr12gb,  $$
we can write the identity (3.10) in the alternative form:
$$
 B=L^2-h^2b^2.
$$

The derivatives
$$
\D fg=\fr qB\lf(\fr14Gb-\fr1{2h}q\rg),  \qquad
\D{J^2}g=-\fr1{h^3}fJ^2 -G \fr qB\lf(\fr14Gb-\fr1{2h}q\rg) J^2,
$$
and
\be
\D Jb=-\fr g2\fr qBJ, \qquad
\D Jq=\fr g2\fr bBJ
\ee
can readily be obtained from (3.12).

{

In many cases it is convenient to use the variables
\be
u_i~:=a_{ij}y^j,
\qquad
v^i~:=y^i+bb^i, \qquad v_m~:=u_m+bb_m=r_{mn}y^n     \equiv a_{mn}v^n.
\ee
The identities
\be
u_iv^i=v_iy^i=q^2,
 \qquad
v_ib^i=v^ib_i=(1-c^2)b,   \qquad
 r_{ij}b^j=(1-c^2)b_i,
\ee
\ses
\be
v_iv^i=q^2+(1-c^2)b^2,  \qquad
\D b{y^i}=b_i, \qquad \D q{y^i}=\fr{v_i}q,
\qquad \D {(q/b)}{y^i}=\fr q{b^2}
e_i
\ee
can readily be verified,
where
\be
e_i=-b_i+\fr b{q^2}v_i
\ee
is the vector showing the property
$
e_iy^i=0.
$
 Particularly, the vector enters the equality
\be
\D{\lf(\fr{F^2}B\rg)} {y^k}=  \fr{gq}B \fr{F^2}B      e_k.
\ee

\ses

By performing required direct calculation,
we find
the representations
\be
y_i=\Bigl(v_i-(b+gq)b_i\Bigr)\fr{F^2}B,
\ee
\ses
\be
g_{ij}=
\biggl[a_{ij}
-\fr g{B}\Bigl ( -   q(b+gq)b_ib_j+q(b_iv_j+b_jv_i)-b\fr{v_iv_j}q\Bigr)\biggr]\fr{F^2}B,
\ee
and
\be
g^{ij}=
\biggl[a^{ij}
+\fr g{B}\Bigl(-bqb^ib^j   +    q(b^iv^j+b^jv^i)  -    (b+gc^2q)\fr{v^iv^j}{\nu}\Bigr)
\biggr]\fr B{F^2}.
\ee
The determinant of the metric tensor  is
\be
\det(g_{ij})= \fr{\nu}q   \lf(\fr{F^2}B\rg)^N\det(a_{ij}),
\ee
where
\be
\nu=q-(1-c^2)gb.
\ee

{

Contracting shows that
$$
g_{ij}b^j=
\Biggl[b_i
-\fr g{B}\lf (
   q(b+gq)c^2b_i  -qc^2v_i
+qb(1-c^2)b_i-\fr bq(1-c^2)bv_i\rg)
\Biggr]\fr{F^2}B,
$$
\ses
or
\be
g_{ij}b^j=
\Biggl[b_i
-\fr g{B}\lf (
   q(b+gc^2q)b_i
                 -\fr 1q(c^2q^2+(1-c^2)b^2)v_i\rg)
\Biggr]\fr{F^2}B,
\ee
\ses
and
$$
g_{ij}v^j=\Biggl[v_i
-(b+gq)b_i
+bb_i
-b\fr g{B}\lf (
   q(b+gc^2q)b_i
 -\fr 1q(c^2q^2+(1-c^2)b^2)v_i
 \rg)
\Biggr]\fr{F^2}B,
$$
\ses
or
\be
g_{ij}v^j=\fr{\nu}q
\Biggl[
(q^2-b^2)v_i
-gq^3b_i
\Biggr]\fr{F^2}{B^2}.
\ee

It is of help  to note that
\be
b(b+gc^2q)=-B+q\nu.
\ee

{

To verify that the tensor (3.21) is reciprocal to (3.20), we are to
demonstrate
that   $g^{nj}g_{ij}-\de^n{}_i =0$. We have
$$
g^{nj}g_{ij}-\de^n{}_i=
-\fr g{B}\Bigl ( -   q(b+gq)b_ib^n+q(b_iv^n+b^nv_i)-b\fr{v_iv^n}q\Bigr)
$$

\ses

$$
+\fr{gq}B(v^n-bb^n)
\Biggl[b_i
-\fr g{B}\lf (
   q(b+gc^2q)b_i
                 -\fr 1q(c^2q^2+(1-c^2)b^2)v_i\rg)
\Biggr]
$$

\ses

$$
+\fr g{B^2q}\Bigl(q\nu b^n-(b+gc^2q)v^n\Bigr)
\lf[
(q^2-b^2)v_i
-gq^3b_i
\rg].
$$
\ses
Simplifying yields
$$
g^{nj}g_{ij}-\de^n{}_i=
-\fr g{B}\Bigl ( q(b_iv^n+b^nv_i)-b\fr{v_iv^n}q\Bigr)
$$

\ses

$$
+\fr{gq}Bv^n
\Biggl[b_i
-\fr g{B}\lf (
   q(b+gc^2q)b_i
                 -\fr 1q(c^2q^2+(1-c^2)b^2)v_i\rg)
\Biggr]
-\fr{gq}Bbb^n
\fr g{B}
\fr 1q(c^2q^2+(1-c^2)b^2)v_i
$$

\ses

$$
+\fr g{B^2q}q\nu b^n
(q^2-b^2)v_i
-\fr g{B^2q}(b+gc^2q)v^n
\lf[
(q^2-b^2)v_i
-gq^3b_i
\rg],
$$
\ses
so the rest is
$$
g^{nj}g_{ij}-\de^n{}_i=
         \fr g{B}b\fr{v_iv^n}q
+\fr{gq}B\fr gBv^n
\fr 1q(c^2q^2+(1-c^2)b^2)v_i
-\fr g{B^2q}(b+gc^2q)v^n
(q^2-b^2)v_i=0.
$$
Thus, the reciprocity is valid.

{

In terms of the set $\{b_i,u_i=a_{ij}y^j \}$, we obtain the alternative representations
\be
y_i=\bigl(u_i - gqb_i\bigr)\fr{F^2}B
\ee
\ses
and
\be
g_{ij}=
\biggl[a_{ij}
+\fr g{B}\Bigl((gq^2-\fr{b(q^2-b^2)}q)b_ib_j  +  \fr bqu_iu_j
-
\fr{ q^2-b^2  }q(b_iu_j+b_ju_i)\Bigr)\biggr]\fr{F^2}B.
\ee

We get also
\be
g^{ij}=
\biggl[a^{ij}+\fr g{\nu}(bb^ib^j+b^iy^j+b^jy^i) -   \fr g{B\nu}(b+gc^2q)y^iy^j
\biggr]\fr B{F^2}.
\ee

With the help of (3.27) we can transform (3.28) to
 the representation
\be
g_{ij}=
\biggl[a_{ij}
+\fr g{q}
\Bigl(
-(b+gq)b_ib_j
+   \fr1B b  \fr B{F^2}   \fr B{F^2}y_iy_j
-  \fr B{F^2}(b_iy_j+b_jy_i)\Bigr)
\biggr]\fr{F^2}B.
\ee

{

From the determinant  value (3.22) we can explicate the vector
\be
 A_i=
 F\D{\ln\lf(\sqrt{|\det(g_{mn})|}\rg)}{y^i},
  \ee
 obtaining
$$
A_i= \fr12F(1-c^2)\fr g{\nu}e_i
+\fr NFy_i
-\fr {NF}{2B}
\Bigl(2v_i-2bb_i-gqb_i-gb\fr1qv_i\Bigr),
$$
\ses
or
\be
A_i= g\fr F{2B}
\fr qX
e_i,
\ee
where
\be
\fr1X=    N+\fr{(1-c^2)B}{q\nu}.
\ee
Also,
$$
A^i=g\fr 1{2F}
\fr qX
\lf(-b^i+\fr b{q^2}v^i+\fr g{q^2\nu}v^i(c^2q^2+(1-c^2)b^2)\rg).
$$
Simplifying  yields
\be
A^i=g\fr 1{2F}
\fr qX
\lf(-b^i  +\fr 1{q\nu}(b+gc^2q)v^i\rg),
\ee
or
alternatively,
\be
A^i=g\fr 1{2F}
\fr 1{X\nu}
\Bigl(-Bb^i  +(b+gc^2q)y^i\Bigr).
\ee
We obtain the contraction
\be
A^hA_h=
  - \fr{g^2}4
\fr1{X^2}\lf(N+1-\fr1X\rg).
\ee

{

Differentiating (3.30) yields
$$
\D{g_{ij}} {y^k}-  \fr2F XA_k g_{ij}
=
\fr gq e_k
                  b_ib_j
\fr{F^2}B
$$

\ses

$$
-
\fr gq\biggl[
  e_kl_il_j
+   2    \Bigl(\fr{1}{F}l_k      -      \fr{1}F XA_k\Bigr)
b  l_il_j
\biggr]
+
\fr gq b    \fr 1{F}(g_{ik}l_j+g_{jk}l_i)
-
\fr gq
(g_{ik}
         b_j
+g_{jk}
       b_i
)
$$

\ses

$$
-
\fr gq
\Bigl(
2 b  \fr 1{F}  l_il_j
-l_jb_i-l_ib_j\Bigr)
2 XA_k
+\fr g{qb}e_k (y_j          b_i +y_i          b_j  )
+\fr g{qb}  y_j
           b_kb_i
+\fr g{qb}
y_i
         b_kb_j.
$$
Applying here the equalities
\be
b_i=bl_i\fr1{F}-\fr {q^2}Be_i
\ee
(see (3.19))
and
\be
b_ib_j=
 b^2l_il_j\fr1{F^2}-bq^2\fr 1B  \fr1{F}  (e_il_j+e_jl_i)+q^4\fr1{B^2}e_ie_j,
\ee
we obtain

{

$$
\D{g_{ij}} {y^k}-  \fr2F XA_k g_{ij}
=
\fr gq e_k
                  b_ib_j
\fr{F^2}B
-
\fr gq
  e_kl_il_j
$$

\ses

$$
+
\fr gq b    \fr 1{F}(h_{ik}l_j+h_{jk}l_i)
-
\fr gq
\Biggl(g_{ik}
         b_j
+g_{jk}
       b_i
\Biggr)
-2\fr gqb\fr1FXA_k   l_il_j
+
2\fr gq
\Bigl(l_jb_i+l_ib_j\Bigr)
 XA_k
$$

\ses

$$
+\fr g{qb}e_k\Biggl(y_j
          b_i
+y_i
         b_j
\Biggr)
+\fr g{qb}  y_j
\Biggl(bl_k\fr1{F}-\fr {q^2}Be_k\Biggr)b_i
+\fr g{qb}
y_i
\Biggl(bl_k\fr1{F}-\fr {q^2}Be_k\Biggr)b_j
$$

\ses

\ses

\ses

$$
=
-
\fr gq
  e_kl_il_j
-2\fr gqb\fr1FXA_k   l_il_j
+
\fr gq e_k
 \Biggl( b^2l_il_j\fr1{F^2}-bq^2\fr 1B  \fr1{F}  (e_il_j+e_jl_i)+q^4\fr1{B^2}e_ie_j\Biggr)
\fr{F^2}B
$$

\ses

$$
+
\fr gq b    \fr 1{F}(h_{ik}l_j+h_{jk}l_i)
-
\fr gq
\Biggl(h_{ik}
         b_j
+h_{jk}
       b_i
\Biggr)
-\fr g{q}
b\fr {1}B
e_k\Biggl(y_j
          b_i
+y_i
         b_j
\Biggr)
$$

\ses

\ses

\ses

$$
=
-
\fr gq
  e_kl_il_j
-2\fr gqb\fr1FXA_k   l_il_j
+
\fr gq e_k
 \Biggl( b^2l_il_j\fr1{F^2}-bq^2\fr 1B  \fr1{F}  (e_il_j+e_jl_i)+q^4\fr1{B^2}e_ie_j\Biggr)
\fr{F^2}B
$$

\ses

\ses

$$
+
\fr gq b    \fr 1{F}(h_{ik}l_j+h_{jk}l_i)
-
\fr gq
\Biggl(h_{ik}
 \Bigl(bl_j\fr1{F}-\fr {q^2}Be_j\Bigr)
+h_{jk}
 \Bigl(bl_i\fr1{F}-\fr {q^2}Be_i\Bigr)
\Biggr)
$$

\ses

$$
-\fr g{q}
b\fr {1}B
e_k\Biggl(y_j
 \Bigl(bl_i\fr1{F}-\fr {q^2}Be_i\Bigr)
+y_i
 \Bigl(bl_j\fr1{F}-\fr {q^2}Be_j\Bigr)
\Biggr),
$$
\ses
so that
\ses\\
\be
\D{g_{ij}} {y^k}-  \fr2F XA_k g_{ij}
=
 gq^3
\fr{F^2}{B^3}   e_ke_ie_j
-
\fr {gq}B
  e_kl_il_j
+
 gq\fr1B
\Bigl(h_{ik}
e_j
+h_{jk}
e_i
\Bigr).
\ee

{

Therefore, the  associated  Cartan tensor
\be
 A_{ijk}~ := \fr F2\D{g_{ij}}{y^k}
 \ee
is representable in the form
\be
A_{ijk}= X
 \Biggl[
A_ih_{jk}  +A_jh_{ik}  +A_kh_{ij}
-\lf(N+1-\fr1X\rg)
\fr1{A_hA^h}A_iA_jA_k
\Biggr].
\ee

\ses

It is useful to verify that
 the contraction
\be
 A_i=g^{jk}A_{ijk}
\ee
when obtained from (3.41) is identical to the representation (3.32).

In various processes of evaluation, it is useful to apply the formulas
\be
\D B{y^k}=\fr{2B}{F^2}y_k      -      \fr{2B}F XA_k,
\qquad
\D{\lf(\fr{F^2}B\rg)} {y^k}=  \fr{2F}B XA_k,
\ee
\ses
and
\be
Fb_n= bl_n  - \fr{2q}{g}X A_n,
 \qquad
F\fr1q v_n=ql_n  +   \fr{B-q^2}b \fr{2}{g} XA_n,
\ee
\ses
together with
\be
\D{\lf( \fr Fq\rg)}{y^n}=    -     \fr{2(B-q^2)}{gbq^2} XA_n,
\qquad
\D{\lf(\fr{q^2}B\rg)} {y^k}=  -      \fr{2q^2  (2b+gq)}{gFB} XA_k,
\ee
which can be verified by
simple straightforward calculation.
It follows that
\be
\D{(XA_k)}{y^n}  =     -  \fr {1}{F}l_k\, XA_n
+\fr{gb}{2F q}  h_{kn}
- \fr{2}{gbqF}(B-q^2)   X^2A_kA_n.
\ee

\ses
\ses

\setcounter{sctn}{4}
         \setcounter{equation}{0}

\nin
{\bf 4.  The space-like  sector of the   space ${\cA\cR}_{g;c}$}

\ses

\ses

Instead of (3.1) and (3.2), we now take

\be
-2<g(x)<2
\ee
and
\be
h(x)=\sqrt{1-\fr14(g(x))^2}.
\ee

In the {\it  space-like  sector}
\be
y\in\cR_g^{+},  \qquad            \cR_g^{+}=\cR_g^{+r}\cup\cR_g^{+l}
\ee
with
\medskip
\be
\cR_g^{+r}=\Bigl(y\in \cR_g^{+r}:~y\in T_xM, \,  b\ge0,\,  q <  -g_-b\Bigr)
\ee
and
\medskip
\be
\cR_g^{+l}=\Bigl(y\in \cR_g^{+l}:~y\in T_xM, \,  b\le0,\,  q<-g_+b\Bigr),
\ee
where
$$
g_+=-\fr12g+h,  \quad  g_-=-\fr12g-h,
$$
 we have
\be
r_{ij}y^iy^j<0
\ee
and can write the variable (2.7) merely as
\be
q=\sqrt{-r_{ij}(x)y^iy^j}.
\ee
The quadratic form (2.9) takes on the form
\be
B =-(q^2+gbq+b^2) < 0,
\ee
which is of the negative discriminant
\be
D_{\{B\}}=-4h^2<0,
\ee
where $h$ is given by (4.2).

{

\ses  \ses

In the  space-like  region (4.3), the squared metric function $F^2(x,y)$
 is given by the formulas
\be
F^2(x,y)=
B(x,y)\,J^2(x,y)
\ee
and
\be
J(x,y)=\e^{-\frac12G(x)f(x,y)}
\ee
with
\be
f=
-\arctan \fr G2+\arctan\fr{L}{hb}
=\arctan\fr {hq}{b+\fr12gq},
\qquad  {\rm if}  \quad b\le 0,
\ee
and
\be
f= - \pi-\arctan\fr G2+\arctan\fr{L}{hb}
=
-\pi+\arctan\fr {hq}{b+\fr12gq},
\qquad  {\rm if}
 \quad b\ge 0,
\ee
where $G=g/h$ and
 \be
 L =q+\fr g2b.
\ee

We have
\be
 f=0,\quad {\rm if} \quad q=0 \quad {\rm and} \quad b<0;
\ee

\ses

\be
f=
-\pi,\quad {\rm if} \quad q=0 \quad {\rm and} \quad b>0.
\ee

Sometimes it is convenient to use also the function
\be
A=b+\fr g2q.
\ee
The identities
\be
L^2+h^2b^2=-B, \qquad h^2q^2+A^2=-B
\ee
are valid.

The derivatives
\be
\D Jb=-\fr g2\fr qBJ, \qquad
\D Jq=\fr g2\fr bBJ
\ee
can readily be obtained from (4.11)-(4.13).

It is useful to compare (4.10)-(4.13) with (3.11) and (3.12).

{

Again, we  use the variables
\be
u_i~:=a_{ij}y^j,
\qquad
v^i~:=y^i+bb^i, \qquad v_m~:=u_m+bb_m=r_{mn}y^n     \equiv a_{mn}v^n,
\ee
and apply the identities
\be
u_iv^i=v_iy^i=-q^2,
 \qquad
v_ib^i=v^ib_i=(1-c^2)b,   \qquad
 r_{ij}b^j=(1-c^2)b_i,
\ee
\ses
\be
v_iv^i=-q^2+(1-c^2)b^2,  \qquad
\D b{y^i}=b_i, \qquad \D q{y^i}=-\fr{v_i}q,
\qquad \D {(q/b)}{y^i}=\fr q{b^2}
e_i,
\ee
where
\be
e_i=-b_i -    \fr b{q^2}v_i
\ee
is the vector showing the property
$
e_iy^i=0.
$
 Particularly, the vector enters the equality
\be
\D{\lf(J^2\rg)} {y^k}=  \fr{gq}B J^2      e_k.
\ee

{

By performing required direct calculation,
we find
the representations
\be
y_i=\Bigl(v_i -   (b+gq)b_i\Bigr)J^2,
\ee
\ses
\be
g_{ij}=
\biggl[a_{ij}
-\fr g{B}\Bigl ( -   q(b+gq)b_ib_j+q(b_iv_j+b_jv_i)  +   b\fr{v_iv_j}q\Bigr)\biggr] J^2,
\ee
and
\be
g^{ij}=
\biggl[a^{ij}
+\fr g{B}\Bigl(-bqb^ib^j   +    q(b^iv^j+b^jv^i)  +    (b+gc^2q)\fr{v^iv^j}{\nu}\Bigr)
\biggr] \fr1{J^2}.
\ee
The determinant of the metric tensor  is
\be
\det(g_{ij})= \fr{\nu}q   \lf(J^2\rg)^N\det(a_{ij}),
\ee
where
\be
\nu=q +  (1-c^2)gb.
\ee

{

Contracting shows that
$$
g_{ij}b^j=
\Biggl[b_i
-\fr g{B}\lf (
   q(b+gq)c^2b_i  -qc^2v_i
+qb(1-c^2)b_i   +     \fr bq(1-c^2)bv_i\rg)
\Biggr]J^2,
$$
\ses
or
\be
g_{ij}b^j=
\Biggl[b_i
-\fr g{B}\lf (
   q(b+gc^2q)b_i
                 -\fr 1q\bigl(c^2q^2  -(1-c^2)b^2\bigr)v_i\rg)
\Biggr]J^2,
\ee
\ses
and
$$
g_{ij}v^j=\Biggl[v_i
-(b+gq)b_i
+bb_i
-b\fr g{B}\lf (
   q(b+gc^2q)b_i
 -\fr 1q\bigl(c^2q^2  -   (1-c^2)b^2\bigr)v_i
 \rg)
\Biggr]J^2,
$$
\ses
or
\be
g_{ij}v^j=      \fr{\nu}q
\Biggl[
-    (q^2+b^2)v_i
+gq^3b_i
\Biggr]\fr1BJ^2.
\ee

It is useful to note that
\be
b(b+gc^2q)=-B-q\nu.
\ee

{

Let us  verify that the tensor (4.27) is reciprocal to (4.26):
$$
g^{nj}g_{ij}-\de^n{}_i=
-\fr g{B}\Bigl ( -   q(b+gq)b_ib^n+q(b_iv^n+b^nv_i) +   b\fr{v_iv^n}q\Bigr)
$$

\ses

$$
+\fr{gq}B(v^n-bb^n)
\Biggl[b_i
-\fr g{B}\lf (
   q(b+gc^2q)b_i
                 -\fr 1q\bigl(c^2q^2-(1-c^2)b^2\bigr)v_i\rg)
\Biggr]
$$

\ses

$$
+\fr g{B^2q}\Bigl(q\nu b^n+(b+gc^2q)v^n\Bigr)
\lf[
-(q^2+b^2)v_i
+gq^3b_i
\rg].
$$
\ses
Simplifying yields
$$
g^{nj}g_{ij}-\de^n{}_i=
-\fr g{B}\Bigl ( q(b_iv^n+b^nv_i)+b\fr{v_iv^n}q\Bigr)
$$

\ses

$$
+\fr{gq}Bv^n
\Biggl[b_i
-\fr g{B}\lf (
   q(b+gc^2q)b_i
                 -\fr 1q\bigl(c^2q^2  -  (1-c^2)b^2\bigr)v_i\rg)
\Biggr]
-\fr{gq}Bbb^n
\fr g{B}
\fr 1q\bigl(c^2q^2  -  (1-c^2)b^2\bigr)v_i
$$

\ses

$$
-\fr g{B^2q}q\nu b^n
(q^2+b^2)v_i
+\fr g{B^2q}(b+gc^2q)v^n
\lf[
-(q^2+b^2)v_i
+gq^3b_i
\rg],
$$
\ses
so the rest is
$$
g^{nj}g_{ij}-\de^n{}_i=
     -    \fr g{B}b\fr{v_iv^n}q
+\fr{gq}B\fr gBv^n
\fr 1q\bigl(c^2q^2  -   (1-c^2)b^2\bigr)v_i
-\fr g{B^2q}(b+gc^2q)v^n
(q^2+b^2)v_i=0.
$$
Thus, the reciprocity is valid.

{

In terms of the set $\{b_i,u_i=a_{ij}y^j \}$, we obtain from  (4.25) and (4.26)
the alternative representations
\be
y_i=\bigl(u_i - gqb_i\bigr)J^2
\ee
\ses
and
\be
g_{ij}=
\biggl[a_{ij}
+\fr g{B}\Bigl((gq^2-\fr{b(q^2+b^2)}q)b_ib_j     -  \fr bqu_iu_j
-
\fr{ q^2+b^2}q(b_iu_j+b_ju_i)\Bigr)\biggr]J^2.
\ee
We get also from (4.27):
\be
g^{ij}=
\biggl[a^{ij}  -  \fr g{\nu}(bb^ib^j+b^iy^j+b^jy^i) +   \fr g{B\nu}(b+gc^2q)y^iy^j
\biggr]\fr 1{J^2}.
\ee
With the help of (4.33) we can transform (4.34) to
\be
g_{ij}=
\biggl[a_{ij}
+\fr g{q}
\Bigl(
(b+gq)b_ib_j
-  \fr1B b   \fr1{J^2} \fr1{J^2}  y_iy_j
+  \fr 1{J^2}(b_iy_j+b_jy_i)\Bigr)
\biggr]J^2.
\ee

{

From the determinant  value (4.28) we can explicate the vector
\be
C_i=
 \D{\ln\lf(\sqrt{|\det(g_{mn})|}\rg)}{y^i},
  \ee
 obtaining
$$
C_i=- \fr12(1-c^2)\fr g{\nu}e_i
+\fr1{F^2} Ny_i
-\fr {N}{2B}
\Bigl(2v_i-2bb_i-gqb_i+gb\fr1qv_i\Bigr),
$$
\ses
or
\be
C_i= g\fr{ 1}  {2B}
\fr qX
e_i,
\ee
where
\be
\fr1X=    N-  \fr{(1-c^2)B}{q\nu}.
\ee
Also,
$$
C^i=g\fr 1{2F^2}
\fr qX
\lf(-b^i-\fr b{q^2}v^i-\fr g{q^2\nu}v^i(c^2q^2   -      (1-c^2)b^2)\rg).
$$
Simplifying  yields
\be
C^i=g\fr 1{2F^2}
\fr qX
\lf(-b^i  -   \fr 1{q\nu}(b+gc^2q)v^i\rg),
\ee
or
alternatively,
\be
C^i=g\fr 1{2F^2}
\fr 1{X\nu}
\Bigl(Bb^i  -    (b+gc^2q)y^i\Bigr).
\ee
We obtain the contraction
\be
F^2C^hC_h=
 \fr{g^2}4
\fr1{X^2}\lf(N+1-\fr1X\rg).
\ee

{

Differentiating (4.36) yields
$$
\D{g_{ij}} {y^k}- 2 X C_k g_{ij}
=
- \fr gq e_k                      b_ib_j    J^2
$$

\ses

$$
+
\fr gq\biggl[
  e_k \fr{1}{F^2} y_iy_j
+   2    \Bigl(\fr{1}{F^2}y_k      -       X C_k\Bigr)
b  \fr{1}{F^2} y_iy_j
\biggr]
-
\fr gq b    \fr 1{F^2}(g_{ik}y_j+g_{jk}y_i)
+
\fr gq
(g_{ik}
         b_j
+g_{jk}
       b_i
)
$$

\ses

$$
+
\fr gq
\Bigl(
2 b  \fr 1{F^2}  y_iy_j
-(y_jb_i+y_ib_j)\Bigr)
2 X C_k
-\fr g{qb}e_k (y_j          b_i +y_i          b_j  )
-\fr g{qb}  y_j
           b_kb_i
-\fr g{qb}
y_i
         b_kb_j.
$$
Applying here the equalities
\be
b_i= by_i\fr1{F^2}+\fr {q^2}Be_i
\ee
(see (4.25))
and
\be
b_ib_j=
 b^2y_iy_j\fr1{F^4}+ bq^2\fr 1B  \fr1{F^2}  (e_iy_j+e_jy_i)+q^4\fr1{B^2}e_ie_j,
\ee
we obtain

{

$$
\D{g_{ij}} {y^k}-  2 X  C_k g_{ij}
=
-  \fr gq e_k                     b_ib_j    J^2
+
\fr gq
  e_k   \fr 1{F^2} y_iy_j
$$

\ses

$$
-
\fr gq b    \fr 1{F^2}(h_{ik}y_j+h_{jk}y_i)
+
\fr gq
\Bigl(g_{ik}
         b_j
+g_{jk}
       b_i
\Bigr)
+2\fr gqb\fr1{F^2}X C_k   y_i  y_j
-
2\fr g{q}
\Bigl(y_jb_i+y_ib_j\Bigr)
 X C_k
$$

\ses

$$
-\fr g{qb}e_k\Biggl(y_j
          b_i
+y_i
         b_j
\Biggr)
-\fr g{qb}  y_j
\Biggl(by_k\fr1{F^2}+\fr {q^2}Be_k\Biggr)b_i
-\fr g{qb}
y_i
\Biggl(by_k\fr1{F^2}+\fr {q^2}Be_k\Biggr)b_j
$$

\ses

\ses

\ses

$$
=
\fr gq
  e_k   \fr 1{F^2} y_iy_j        +  2\fr gqb\fr1{F^2}X C_k   y_iy_j
-
\fr gq e_k
 \Biggl(
 b^2y_iy_j\fr1{F^4}+ bq^2\fr 1B  \fr1{F^2}  (e_iy_j+e_jy_i)+q^4\fr1{B^2}e_ie_j  \Biggr)
J^2
$$

\ses

$$
-
 \fr gq b    \fr 1{F^2}(h_{ik}y_j+h_{jk}y_i)
+
\fr gq
\Biggl(h_{ik}
         b_j
+h_{jk}
       b_i
\Biggr)
+\fr g{q}
b\fr {1}B
e_k\Biggl(y_j
          b_i
+y_i
         b_j
\Biggr)
$$

\ses

\ses

                  \ses

$$
=
\fr gq
  e_k   \fr 1{F^2} y_iy_j        +  2\fr gqb\fr1{F^2}X C_k   y_iy_j
-
\fr gq e_k
 \Biggl(
 b^2y_iy_j\fr1{F^4}+ bq^2\fr 1B  \fr1{F^2}  (e_iy_j+e_jy_i)+q^4\fr1{B^2}e_ie_j \Biggr)
J^2
$$

\ses

\ses

$$
-
 \fr gq b    \fr 1{F^2}(h_{ik}y_j+h_{jk}y_i)
+
\fr gq
\Biggl(h_{ik}
 \Bigl(by_j\fr1{F^2}+\fr {q^2}Be_j\Bigr)
+h_{jk}
 \Bigl(by_i\fr1{F^2}+\fr {q^2}Be_i\Bigr)
\Biggr)
$$

\ses

$$
+\fr g{q}
b\fr {1}B
e_k\Biggl(y_j
 \Bigl(by_i\fr1{F^2}+\fr {q^2}Be_i\Bigr)
+y_i
 \Bigl(by_j\fr1{F^2}+\fr {q^2}Be_j\Bigr)
\Biggr),
$$
\ses
so that
\ses\\
\be
\D{g_{ij}} {y^k}- 2 X C_k g_{ij}
=
- gq^3
\fr{F^2}{B^3}   e_ke_ie_j
-
\fr {gq}B \fr1{F^2}
  e_ky_iy_j
+
 gq\fr1B
\Bigl(h_{ik}
e_j
+h_{jk}
e_i
\Bigr).
\ee
This shows that the  associated  Cartan tensor
is representable in the form (2.48).
We have applied the equalities

{

\be
\D B{y^k}=\fr{2B}{F^2}y_k      -      2BX C_k,
\qquad
\D{\lf(J^2\rg)} {y^k}=  \fr{2}B X F^2 C_k.
\ee

{

\ses

\ses

\setcounter{sctn}{5}
 \setcounter{equation}{0}

\nin {\bf 5.  Fixing the tangent space
}

\ses  \ses

Let us introduce the {\it orthonormal frame}
  $h^p_i(x)$ of the input pseudo-Riemannian metric tensor
$a_{ij}(x)$:
\be
a_{ij}=e_{pq}h^p_ih^q_j,
\ee
where $\{e_{pq}\}$ is the pseudo-Euclidean diagonal:
\be
e_{pq}=\text{diagonal}(+ - ...-);
\ee
the indices
$p,q,...$ will be specified on the range $0,1,...,N-1$.
 Denote by $h_p^i$ the reciprocal frame,
 so that
 \be
 h_p^jh^p_i=\de^j_i.
 \ee
At any fixed point $x$, we can represent the tangent vectors $y$ by their frame-components:
\be
R^p=h^p_iy^i,
\ee
and consider the respective  metric tensor components
\be
g_{pq}=h_p^ih_q^jg_{ij}
\ee
of the Finslerian metric tensor $g_{ij}$. We obtain the decomposition
\be
R^p=\{R^a,R^{N-1}\},
\ee
often denoting
\be
R^{N-1}=z.
\ee
The indices
$a,b,...$ will be specified on the range $1,...,N-1$.

It is convenient to specify the frame such that the $(N-1)$-th component $h^{N-1}_i(x)$
is collinear to the input vector field $b_i(x)$.
Under these conditions  the input 1-form $b$ reads
\be
b=cz
\ee
and we have
\be
b^p=\{0,0,...,-c\}, \qquad    b_p=\{0,0,...,c\},
\ee
together with
\be
 q^2=\epsilon\Bigl[e_{ad}R^aR^d - (1-c^2)(R^{N-1})^2\Bigr].
\ee
Notice the inequality
\be
\epsilon\Bigl[e_{ad}R^aR^d - (1-c^2)(R^{N-1})^2\Bigr]>0,
\ee
so that $q$ is a nowhere vanishing variable in both the time-like and space-like sectors.

From (2.37) it follows that the covariant components $R_p=h_p^iy_i$ are equal to
\be
R_a=
e_{ab}R^b
J^2,
\qquad
R_{N-1}=-\fr1c(b+gc^2q)
J^2.
\ee

In the four-dimensional relativistic space-time,
\be
N=4, \qquad  R^p=\{R^0,R^1,R^2,R^3\},
\ee
 it is convenient to relabel the coordinates $R^p$ as follows:
\be
R^0=t, \quad  R^1=x, \quad R^2=y, \quad R^3=z.
\ee
We get
\be
q=\sqrt{\epsilon\Bigl(t^2-x^2-y^2-(1-c^2)z^2\Bigr)}
\ee
and
\be
B=\epsilon q^2-c^2z^2   -    gcz q  \equiv t^2 -x^2 - y^2  - z^2    -    gcz q.
\ee
In terms of such coordinates, the metric tensor components
$g_{pq}$ are obtained
 from  the tensorial components $g_{ij}$ given in (2.38).
The result reads
\be
g_{00}=  \Bigl(   1+\epsilon\fr {gc}{Bq}z t^2   \Bigr)  J^2, \qquad
g_{01}=    -\epsilon\fr {gc}{Bq}zx t    J^2,
\qquad
g_{02}=    -\epsilon\fr {gc}{Bq}zy t    J^2,
\ee

\ses

\be
g_{03}=  - \fr {gc}{Bq}E      t J^2,  \qquad
  g_{11}=\Bigl(  - 1 +\epsilon \fr {gc}{Bq}zx^2   \Bigr)  J^2,   \qquad
g_{22}=\Bigl( -  1 + \epsilon\fr {gc}{Bq} zy^2   \Bigr)  J^2,  \qquad
\ee

\ses

\be
g_{12}= \epsilon\fr {gc}{Bq}zxy    J^2, \qquad
g_{13}=    \fr {gc}{Bq}
E  x  J^2, \qquad
g_{23}=\fr {gc}{Bq} E  y  J^2,
\ee

\ses

\be
g_{33}=\Biggl[
-1+\fr {gc} {Bq}  \Bigl(\bigl(gq^3-b(q^2-\epsilon b^2)\bigr)c  + \epsilon  z^3
+
2( q^2 -\epsilon b^2)z\Bigr)   \Biggr]  J^2.
\ee
Here, $E=\epsilon z^2    + (q^2-\epsilon b^2)$;
$\epsilon = 1$ in the    future-time-like sector    $   y\in\cT_g^{+}$
and $\epsilon = -1$   in the   space-like  sector    $ y\in\cR_g^{+}.$
From (2.39) we can obtain the contravariant components
\be
g^{00}=  \Bigl(   1-\epsilon\fr {g}{B\nu} {\beta} t^2   \Bigr)  \fr1{J^2}, \qquad
g^{01}=    -\epsilon\fr {g}{B\nu} {\beta}x t    \fr1{J^2},
\qquad
g^{02}=    -\epsilon\fr {g}{B\nu} {\beta}y t    \fr1{J^2},
\ee

\ses

\be
g^{03}=  - \epsilon\fr {g}{B\nu}\lf( Bc    + {\beta}z   \rg)t  \fr1{J^2},
\quad
g^{11}=\Bigl(  - 1-\epsilon\fr {g}{B\nu} {\beta} t^2   \Bigr)  \fr1{J^2},
   \quad
g^{22}= \Bigl( -  1-\epsilon\fr {g}{B\nu} {\beta} t^2   \Bigr)  \fr1{J^2},
\ee

\ses

\be
g^{12}=- \epsilon\fr {g}{B\nu} {\beta}yx\fr1{J^2},
\quad
g^{13}=    - \epsilon\fr {g}{B\nu}\lf( Bc    + {\beta}z   \rg)x  \fr1{J^2},
 \quad
g^{23}= - \epsilon\fr {g}{B\nu}\lf( Bc    + {\beta}z   \rg)y  \fr1{J^2},
\ee

\ses

\be
g^{33}=\Biggl[-  1  +  \epsilon\fr {g}{\nu}(c^2-2)cz
-\epsilon\fr {g}{B\nu} {\beta} z^2
\Biggr]  \fr1{J^2},
\ee
where $\beta=b+gc^2q$.

At an arbitrary dimension $N\ge2$,
from (2.38) we obtain the respective components
\be
g_{ab}=
\lf[
e_{ab} + \epsilon g\fr{  e_{ad}R^d   e_{be}R^e   b   }{Bq}\rg]
J^2,  \quad
g_{N-1,a}=\fr{g}{Bq}\biggl(-\epsilon bz-(q^2-\epsilon b^2)c\biggr)
e_{ab}R^b
J^2,
\ee

\ses

\ses

\be
 g_{N-1,N-1}=\Biggl[-1+\fr g{Bq}\biggl(  \bigl(gq^3-b (q^2-\epsilon b^2)\bigr)c^2
+\epsilon bz^2
+2(q^2-\epsilon b^2)b     \biggr)
\Biggr]
J^2.
\ee
\ses
We may alternatively obtain them through the rule
$
g_{pq}
=\partial R_q/\partial R^p.
$
From (2.39) it follows that
\be
g^{ab}=
\lf[
e^{ab} - \epsilon \fr{g}{B\nu}(b+gc^2q)R^aR^b
\rg]\fr1{J^2}, \quad
g^{N-1,a}= -
\epsilon\fr{g}{B\nu}\Bigl[
cB+(b+gc^2q)z
\Bigr]R^a
\fr1{J^2},
\ee

\ses

\be
 g^{N-1,N-1}=\Biggl[-1+ \epsilon\fr{g}{B\nu}\Bigl[
 (bc^2-2b)B
-(b+gc^2q)z^2
\Bigr]
\Biggr]
\fr1{J^2}.
\ee

{

\ses

\ses

\setcounter{sctn}{6} \setcounter{equation}{0}

\nin
  {\bf 6.  Hamiltonian function}

\ses

\ses

It proves  possible to obtain explicitly the {\it Hamiltonian function} $H$
of the  $ {\mathbf\cA\cR}$-space.
  To this end
 we introduce the co-form
$
\hat b~:= y_ib^i.
$
From (2.29) (the formula (2.37) can also be used) we obtain
\be
\hat b=(b+qgc^2)   J^2.
\ee
Using the tensor
$
r^{ij}=a^{ij}+b^ib^j
$
(see (2.5)),  it is natural to introduce the co-counterpart
\be
\hat \ga=y_iy_jr^{ij}
\ee
to the scalar (2.6).  Inserting here (2.37) yields
\be
\hat \ga= \Bigl[\ga  -  2(1-c^2)gbq - c^2(1-c^2)g^2q^2\Bigr]J^4
\ee
(the formulas (2.25)-(2.26)  have been used).

It is appropriate to use  the co-counterpart
\be
\hat q=\sqrt{|\hat\ga|}
\ee
to (2.7)
and rewrite (6.3) as follows:
\be
\hat q^2= \Bigl[q^2  -2(1-c^2)\epsilon gbq -  c^2(1-c^2)\epsilon g^2q^2\Bigr]J^4.
\ee

On resolving the equation set (6.1)-(6.5) to get the functions $b=b(\hat b, \hat q)$
and
 $q=q(\hat b, \hat q)$, we can insert  the functions in the $F^2$ to obtain
 the squared Hamiltonian function $H^2$
 through
the equality
\be
H^2=F^2
\ee
(we apply the method  presented in Section 7.1 of [1]).
This  leads to the Hamiltonian function
 representation
\be
H^2(x,\hat  y) =\hat \Phi \Bigl(g(x), b^i(x), a^{ij}(x),\hat  y\Bigr),
\ee
which is the co-counterpart to (1.1).

The knowledge of $H^2$ makes it possible to formulate
the  $ {\mathbf\cA\cR}$-Hamilton-Jacobi equation
\be
H^2\lf(x, \D S {x^i}\rg)=m^2,
\ee
where $S=S(x)$ is the characteristic function and $m$ is the rest mass of the particle.

When $c=1$,  from (6.1) and (6.5) it ensues that
\be
\hat b=(b+gq)J^2, \qquad \hat q=qJ^2,
\ee
which entails
\be
\fr{\hat q}{\hat b}=\fr q{b+gq  }, \qquad
\fr{ q}{ b}=\fr{\hat q}{\hat b  - g \hat q  }.
\ee
Introducing the quadratic form
\be
\hat B= \epsilon\hat q^2+g \hat  b\hat q - \hat b^2,
\ee
we get  the equality
\be
\hat B=BJ^4,
\ee
where $B$ is the initial quadratic form (2.9).
When the co-vectors $\hat y$ are
relatable to the  vectors $y$ of the future-time-like sector, so that $\epsilon =1$,
from (6.9) it follows that
the   function $H^2(x,\hat y)$
 is obtained from (3.11)-(3.12) to read
\be
H^2(x,\hat y) =  \hat B(x,\hat y)\hat J^2(x,\hat y)
\equiv
(\hat b+g_+\hat q)^{-G_-}  (-g_-\hat q-\hat b)^{G_+},
 \quad\epsilon =1,
\ee
\ses
and
\be
\hat J(x,\hat y)=
\lf(
\fr{\hat b+g_+\hat q}{-g_-\hat q-\hat b}
\rg)^{G/4} \equiv \fr1{J(x,y)},  \quad\epsilon =1,
\ee
with
$G=g/h$, \, $G_+=g_+/h$, and $G_-=g_-/h$;
the quantities
  $g_+$ and $  g_-$
  have been indicated below (4.5).
The quadratic form
\be
\hat B=(\hat b+g_+\hat q)  (-g_-\hat q-\hat b), \quad \epsilon =1,
\ee
substitutes now the initial $B$.

\ses

In the  space-like  sector,  from (4.10)-(4.13) and (6.9)-(6.12) we obtain
\be
H^2(x,\hat y)=
\hat B(x,\hat y)\,\hat J^2(x,\hat y),  \quad\epsilon =-1,
\ee
and
\be
\hat J(x,\hat y)=\e^{\frac12G(x)\hat f(x,\hat y)} \equiv \fr1{J(x,y)}, \quad\epsilon =-1,
\ee
with
\be
\hat f=
\arctan \fr G2+\arctan\fr{\hat L}{h\hat b},
\qquad  {\rm if}  \quad \hat b\ge 0,
\ee
and
\be
\hat f= \pi+\arctan\fr G2+\arctan\fr{\hat L}{h\hat b},
\qquad  {\rm if}
 \quad \hat b\le 0,
\ee
where
 \be
\hat  L =\hat q-\fr g2\hat b.
\ee

We have taken into account  the identities
\be
\arctan\fr{\hat L}{h\hat b}
-\arctan\fr{ L}{h  b}
=
-\arctan\fr{gh}{1-\fr{g^2}2}
\ee
and
\be
\arctan\fr{ G}2
-\arctan\fr{gh}{1-\fr{g^2}2}
=-
\arctan\fr{ G}2.
\ee

{

\ses

\ses

\setcounter{equation}{0}

\nin
{ \bf Appendix A:  Applying the
${\mathbf\cU\cA\cR}$-space  coordinates
}

\ses\ses

Below the coordinates $z^p$ defined in (2.79) are used.

We start with  examining the {\it future-time-like sector}.
 We have $F>0$ and can take $z^0=F$.
We introduce the representations
\ses
\be
R^0=F\cosh\eta \Ch\chi,
 \ee
\ses
\be
R^1=F\sinh\eta  \Ch\chi\cos\phi,
\qquad
R^2=F\sinh\eta \Ch\chi \sin\phi,
\qquad
R^3=F \Sh\chi ,
\ee
obtaining
\be
q=F\Ch\chi,
\ee
with
\be
\Ch\chi=\fr1{Jh} \cosh f, \qquad
\Sh\chi=
\fr1J\lf(\sinh f-\fr G2\cosh f\rg).
\ee
We also introduce the function
\be
\Sh^*\chi=\fr1J\lf(\sinh f + \fr G2\cosh f\rg).
\ee
From this it follows
\be
\Sh'\chi=\Ch\chi,  \qquad    \Ch'\chi=\Sh^*\chi,
\ee
where the prime $'$ means differentiation with respect to the angle $\chi$.
The angle is taken as follows:
\be
\chi =\fr1hf, \qquad
J=e^{-\frac12g\chi}
\ee
(see the formulas below (3.12)).

The identities
\be
\Sh^*\chi=\Sh\chi +  g  \Ch\chi
\ee
and
\be
\Sh^*\chi \Sh^*\chi   -g\Sh^*\chi \Ch\chi   -      \Ch^2\chi
=-\fr1{J^2}
\ee
are valid.

{

Evaluating the partial derivatives yields
\be
\D{R^p}{z^0}=\fr1FR^p,
\ee

\ses

\be
\D{R^0}{z^1}=F\sinh\eta \Ch\chi,\qquad \D{R^0}{z^2}=0,
\qquad
\D{R^0}{z^3}=F\cosh\eta \Sh^*\chi,
\ee

\ses

\be
\D{R^1}{z^1}=F\cosh\eta \Ch\chi \cos\phi,\qquad \D{R^1}{z^2}= - F\sinh\eta\Ch\chi  \sin\phi,
\ee

\ses

\be
\D{R^1}{z^3}=F\sinh\eta \Sh^*\chi \cos\phi,
\qquad
\D{R^2}{z^3}=F\sinh\eta \Sh^*\chi \sin\phi         ,
\ee

\ses

\ses

\be
\D{R^2}{z^1}=F\cosh\eta \Ch\chi \sin\phi,\qquad \D{R^2}{z^2}=  F\sinh\eta\Ch\chi  \cos\phi,
\ee

\ses

\be
\D{R^3}{z^1}=\D{R^3}{z^2}=0,
\qquad
\D{R^3}{z^3}=F\Ch\chi.
\ee

\ses

The derivatives fulfill the identities
\be
R^1\D{R^1}{z^2}
+
R^2\D{R^2}{z^2}
\equiv 0,
 \qquad
\D{R^1}{z^2}
\D{R^1}{z^3}
+
\D{R^2}{z^2}
\D{R^2}{z^3}
\equiv 0.
\ee

{

Let us apply the transformation (A.1)-(A.3) to  the ${\mathbf\cU\cA\cR}$-angle (2.71).
This yields
\be
 \al_{\{x\}}(y_1,y_2)
  = \fr1h\arccosh \tau_{12},
 \ee
where
\be
\tau_{12}=   \cosh f_1  \cosh f_2 Z_{12}-\sinh f_1  \sinh f_2
\ee
with
\be
Z_{12}= \cosh \eta_1 \cosh\eta_2 -\sinh\eta_1 \sinh\eta_2 \cos(\phi_2-\phi_1) .
\ee
We can write also
\be
\tau_{12}=
\cosh(f_2-f_1)+\cosh f_1\cosh f_2
\Bigl(
\cosh(\eta_2-\eta_1) -1 +Y\sinh \eta_1 \sinh\eta_2
\Bigr),
\ee
where $Y=1-\cos(\phi_2-\phi_1)$.

Let us use
 the Finslerian metric tensor components  $g_{pq}(g;R)$
indicated in (5.17)-(5.20) to obtain the tensor
$A_{rs}(g;z)$ introduced by the transformation (2.82).
We get
the diagonal tensor:
$
A_{01}=A_{02}=A_{03}=A_{12}= A_{13}= A_{23}=0
$
and the components
\be
A_{00}=1,\qquad    A_{33}=-(z^0)^2,
\ee
\ses
\be
A_{11}=-(z^0)^2 \fr1{h^2}\cosh^2f,\qquad
A_{22}=-(z^0)^2\sinh^2 \eta \fr1{h^2}\cosh^2f,
\ee
so that
the respective  squared linear element  is
\be
(ds)^2=(dz^0)^2-(z^0)^2\Bigl[
(d\chi)^2+\fr1{h^2}\cosh^2 f
\bigl(\sinh^2\eta\,(d\phi)^2+(d\eta)^2\bigr)
\Bigr].
\ee

Let us verify  the involved coefficients.
We have
$$
A_{11}\!= \!\lf(
\D{R^0}{z^1}g_{00}
+2\D{R^1}{z^1}g_{10}
+2\D{R^2}{z^1}g_{20}
\rg)
\!
\D{R^0}{z^1}
+\D{R^1}{z^1}\D{R^1}{z^1}g_{11}
+2\D{R^1}{z^1}\D{R^2}{z^1}g_{12}
+\D{R^2}{z^1}\D{R^2}{z^1}g_{22},
$$
\ses
so  that
$$
\fr1{J^2F^2}A_{11}=
\sinh^2\eta\Ch^2\chi
\Bigl(   1+\fr {g}{Bq}z t^2   \Bigr)
+  \cosh^2\eta\Ch^2\chi \cos^2\phi \Bigl(  - 1 + \fr {g}{Bq}zx^2   \Bigr)
 $$

\ses

$$
-2\fr {g}{Bq}zx t \sinh\eta\Ch\chi   \cosh\eta\Ch\chi \cos\phi
  -2\fr {g}{Bq}zy t   \sinh\eta\Ch\chi   \cosh\eta\Ch\chi \sin\phi
  $$

\ses

$$
+2  \cosh^2\eta\Ch^2\chi \cos\phi    \sin\phi
\fr {g}{Bq}zxy
+
 \cosh^2\eta\Ch^2\chi \sin^2\phi
\Bigl( -  1 + \fr {g}{Bq} zy^2   \Bigr)
 $$

\ses

\ses

$$
=
\sinh^2\eta\Ch^2\chi
\Bigl(   1+\fr {g}{Bq}z t^2   \Bigr)
-2\fr {g}{Bq}z t F\sinh^2\eta\Ch^2\chi   \cosh\eta\Ch\chi
  $$

\ses

$$
+  \cosh^2\eta\Ch^2\chi
\Bigl(  - 1 + \fr {g}{Bq}zF^2\sinh^2\eta\Ch^2\chi   \Bigr).
$$
By reducing we arrive at
$$
\fr1{J^2F^2}A_{11}=  -\Ch^2\chi.
$$
Therefore,
\be
A_{11}= -  F^2\fr1{h^2}\cosh^2f.
\ee

{

The similar chain:
$$
A_{22}=
\D{R^1}{z^2}\D{R^1}{z^2}g_{11}
+2\D{R^1}{z^2}\D{R^2}{z^2}g_{12}
+\D{R^2}{z^2}\D{R^2}{z^2}g_{22}
$$
and
\ses
$$
\fr1{J^2F^2}A_{22}=
  \sinh^2\eta\Ch^2\chi \sin^2\phi
\Bigl(  - 1 + \fr {g}{Bq}zx^2   \Bigr)
$$

\ses

$$
-2  \sinh^2\eta\Ch^2\chi \cos\phi    \sin\phi
\fr {g}{Bq}zxy
+
 \sinh^2\eta\Ch^2\chi \cos^2\phi
\Bigl( -  1 + \fr {g}{Bq} zy^2   \Bigr)
 $$
leads to
\ses\\
$$
\fr1{J^2F^2}A_{22}=  -\sinh^2\eta\Ch^2\chi,
$$
or
\be
A_{22}= -  F^2\sinh^2 \eta \fr1{h^2}\cosh^2f.
\ee

Next, we consider the component
$$
A_{33}= \lf(\D{R^0}{z^3}g_{00}
+2\D{R^1}{z^3}g_{10}
+2\D{R^2}{z^3}g_{20}
+2\D{R^3}{z^3}g_{30}
\rg)\D{R^0}{z^3}
+\D{R^1}{z^3}\D{R^1}{z^3}g_{11}
+2\D{R^1}{z^3}\D{R^2}{z^3}g_{12}
$$

\ses

$$
+2\D{R^1}{z^3}\D{R^3}{z^3}g_{13}
+\D{R^2}{z^3}\D{R^2}{z^3}g_{22}
+2\D{R^2}{z^3}\D{R^3}{z^3}g_{23}
+\D{R^3}{z^3}\D{R^3}{z^3}g_{33},
$$
\ses
getting
$$
\fr1{J^2F^2}A_{33}=
\cosh^2\eta\Sh^*\chi \Sh^*\chi
\Bigl(   1+\fr {g}{Bq}z t^2   \Bigr)
 $$

\ses

$$
-2\fr {g}{Bq}zx t \cosh\eta\Sh^*\chi \Sh^*\chi   \sinh\eta \cos\phi
  -2\fr {g}{Bq}zy t \cosh\eta\Sh^*\chi \Sh^*\chi   \sinh\eta \sin\phi
  $$

\ses

$$
-2 \cosh\eta\Sh^*\chi \Ch\chi
\fr {g}{Bq}\Bigl(   z^2  +  ( q^2 - b^2)\Bigr)  t
+  \sinh^2\eta\Sh^*\chi \Sh^*\chi   \cos^2\phi
\Bigl(  - 1 + \fr {g}{Bq}zx^2   \Bigr)
$$

               \ses

$$
+2 \sinh^2\eta\Sh^*\chi \Sh^*\chi   \cos\phi    \sin\phi
\fr {gzxy}{Bq}
+2 \sinh\eta\Ch\chi \Sh^*\chi   \cos\phi
   \fr {g  x}{Bq}\Bigl(   z^2  +  ( q^2 - b^2)\Bigr)
$$

\ses

$$
+
 \sinh^2\eta   \Sh^*\chi  \Sh^*\chi  \sin^2\phi
\Bigl( -  1 + \fr {gzy^2  }{Bq}  \Bigr)
+2 \sinh\eta\Sh^*\chi \Ch\chi      \sin\phi
\fr {gy}{Bq} \Bigl(   z^2  + ( q^2 - b^2) \Bigr)
$$

$$
+
    \Ch^2\chi
       \Biggl[
-1+\fr {g} {Bq}  \Bigl(\bigl(gq^3-b(q^2- b^2)\bigr)  +  z^3
+
2( q^2 - b^2)z\Bigr)   \Biggr] .
 $$

\ses

{

\nin
Due simplifying yields
$$
\fr1{J^2F^2}A_{33}=
\cosh^2\eta\Sh^*\chi \Sh^*\chi
\Bigl(   1+\fr {g}{Bq}z t^2   \Bigr)
$$

\ses

$$
-2\fr {g}{Bq}F z t \cosh\eta \sinh^2\eta \Sh^*\chi  \Sh^*\chi  \Ch\chi
-2 \cosh\eta\Sh^*\chi  \Ch\chi
\fr {g}{Bq} q^2  t
  $$

\ses

$$
-  \sinh^2\eta\Sh^*\chi \Sh^*\chi
+  \sinh^2\eta\Sh^*\chi \Sh^*\chi
 \fr {g}{Bq}zF^2\sinh^2 \Ch\chi \Ch\chi
$$

\ses

$$
+2 \sinh\eta\Ch\chi \Sh^*\chi
   \fr {g}{Bq}q^2  F\sinh\eta \Ch\chi
+
    \Ch^2\chi
       \Biggl[
-1+\fr {g} {Bq} q^2(b+gq)
   \Biggr]
 $$

\ses

$$
=
\Sh^*\chi \Sh^*\chi
\Bigl(   1+\fr {gz}{Bq}F^2 \Ch\chi \Ch\chi  \Bigr)
-2 \Sh^*\chi  \Ch\chi
\fr {gq}{B}   F \Ch\chi
-
    \Ch^2\chi
       \Bigl[
1-\fr {gq} {B} (b+gq)
 \Bigr].
 $$
\ses
In this way we obtain simply
$$
\fr1{J^2F^2}A_{33}=
\Sh^*\chi \Sh^*\chi
+\Sh^*\chi \Sh^*\chi
\fr{gbq}B
-2 \Sh^*\chi  \Ch\chi
\fr {gq^2}{B}
-      \Ch^2\chi
+  \Ch\chi         \fr {gq^2} {B} \fr1F (b+gq)
 $$

\ses

$$
=
\Sh^*\chi \Sh^*\chi
-g\Sh^*\chi \Ch\chi
-      \Ch^2\chi
=-\fr 1{J^2}
$$
(the identity (A.9) has been used),
so that
\be
A_{33}= -F^2.
\ee
Thus, the representations (A.21) and (A.22) are correct.

\ses

By making the substitution
\be
z^0=\e^{\si},
\ee
from (A.23) we get
\be
(ds)^2=\e^{2\si}
\Biggl[
(d\si)^2-
(d\chi)^2-\fr1{h^2}\cosh^2 f
\bigl(\sinh^2\eta\,(d\phi)^2+(d\eta)^2\bigr)
\Biggr].
\ee
\ses
Using here new coordinates
\be
\rho=\e^{h\si}\cosh(h\chi),
\qquad
\tau=\e^{h\si}\sinh(h\chi)
\ee
just shows the property that the space under study
is {\it conformally flat}:
\be
(ds)^2= \Ka  \Bigl((ds)^2\Bigr){\Bigl|\Bigr.}_{\text{pseudo-Euclidean}},
\ee
where
\be
\Bigl((ds)^2\Bigr){\Bigl|\Bigr.}_{\text{pseudo-Euclidean}}=
(d\rho)^2-(d\tau)^2
-\rho^2
\Bigl(\sinh^2\eta\,(d\phi)^2+(d\eta)^2\Bigr)
\ee
and
\be
\Ka=\fr1{h^2}(\rho^2-\tau^2)^{(1-h)/h}.
\ee
If we apply here (A.27) and remind that $z^0=F$, we observe that the conformal multiplier
$\Ka$
can be  expressed through the metric function, namely we get
\be
\Ka=\fr1{h^2}\lf(F^2\rg)^{1-h}.
\ee

  Turning to  the {\it space-like sector}, so that
$
\epsilon=-1,
$
we introduce the function
\be
K=\sqrt{|F^2|}>0, \quad \text{so that} \quad z^0=K,
\ee
and set forth
the representations
\be
R^0=K\sinh\eta \Sin\chi,
 \ee
\ses
\be
R^1=K\cosh\eta  \Sin\chi\cos\phi,
\quad
R^2=K\cosh\eta \Sin\chi \sin\phi,
\quad
R^3=K \Cos\chi,
\ee
which entails
\be
q=K\Sin\chi,
\ee
with
\be
\Sin\chi=\fr1{Jh} \sin f, \qquad
\Cos\chi=
\fr1J\lf(\cos f-\fr G2\sin f\rg).
\ee
We need also the function
\be
\Cos^*\chi=\fr1J\lf(\cos f + \fr G2\sin f\rg).
\ee
The derivatives
\be
\Sin'\chi=\Cos^*\chi,  \qquad    \Cos'\chi=-\Sin\chi
\ee
can readily be verified.
Again, we take the angle
$\chi $
according to (A.7)
(notice (4.11)).
There arise the identities
\be
\Cos^*\chi=\Cos\chi +  g  \Sin\chi,\qquad
\Sin^2\chi    +g\Sin\chi \Cos\chi     +    \Cos^2\chi
=\fr1{J^2},
\ee
and also
\be
-\Sin^2\chi    +g\Sin\chi \Cos^*\chi     -\Cos^*\chi \Cos^*\chi=-\fr 1{J^2}.
\ee

{

From (A.35)-(A.36)  we calculate  the partial derivatives, obtaining
\be
\D{R^p}{z^0}=\fr1KR^p,
\qquad
\D{R^3}{z^1}=\D{R^3}{z^2}=0,
\qquad
\D{R^3}{z^3}= - K\Sin\chi,
\ee

\ses

\be
\D{R^0}{z^1}=K\cosh\eta \Sin\chi,\qquad \D{R^0}{z^2}=0,
\qquad
\D{R^0}{z^3}=K\sinh\eta \Cos^*\chi,
\ee

\ses

\be
\D{R^1}{z^1}=K\sinh\eta \Sin\chi \cos\phi,\qquad \D{R^1}{z^2}= - K\cosh\eta\Sin\chi  \sin\phi,
\ee

\ses

\be
\D{R^1}{z^3}=K\cosh\eta \Cos^*\chi \cos\phi,
\qquad
\D{R^2}{z^3}=K\cosh\eta \Cos^*\chi \sin\phi         ,
\ee

\ses

\ses

\be
\D{R^2}{z^1}=K\sinh\eta \Sin\chi \sin\phi,\qquad \D{R^2}{z^2}=  K\cosh\eta\Sin\chi  \cos\phi.
\ee

\ses

Using the substitution (A.35)-(A.37) in   the ${\mathbf\cU\cA\cR}$-angle (2.71) yields
\be
 \al_{\{x\}}(y_1,y_2)
  = \fr1h\arccosh \tau_{12},
 \ee
where
\be
\tau_{12}=   \sin f_1  \sin f_2 Z_{12}+\cos f_1  \cos f_2
\ee
and
\be
Z_{12}=\cosh\eta_1 \cosh\eta_2 \cos(\phi_2-\phi_1)-\sinh \eta_1 \sinh\eta_2.
\ee
Whence we have
\be
\tau_{12}=
\cos(f_2-f_1)+\sin f_1\sin f_2
\Bigl(
\cosh(\eta_2-\eta_1) -1
-  P \cosh \eta_1 \cosh\eta_2
\Bigr),
\ee
where $P=1-\cos(\phi_2-\phi_1)$.

Again, we  evaluate the  transform
$
A_{rs}(g;z)
$
according to (2.82),
now taken the components
(5.17)-(5.20) with the value $\epsilon =-1$.
We get the diagonal tensor: $
A_{01}=A_{02}=A_{03}=A_{12}= A_{13}= A_{23}=0
$
and
\be
A_{00}=-1,\qquad    A_{33}=-(z^0)^2,
\ee
\ses
\be
A_{11}=(z^0)^2 \fr1{h^2}\sin^2f,\qquad
A_{22}=-(z^0)^2\cosh^2 \eta \fr1{h^2}\sin^2f,
\ee
so that the object  $ds^2= - A_{rs}dz^rdz^s$ takes on the form
\be
(ds)^2=-(z^0)^2\Bigl[
\fr1{h^2}\sin^2 f
\bigl((d\eta)^2-\cosh^2\eta\,(d\phi)^2\bigr)
-(d\chi)^2
\Bigr]
+(dz^0)^2.
\ee

\ses

Let us  check the validity of the  key components.
We get
$$
\fr1{J^2K^2}A_{11}=
\cosh^2\eta\Sin^2\chi
\Bigl(   1-\fr {g}{Bq}z t^2   \Bigr)
+  \sinh^2\eta\Sin^2\chi \cos^2\phi
\Bigl(  - 1 - \fr {g}{Bq}zx^2   \Bigr)
 $$

\ses

$$
+2\fr {g}{Bq}zx t \sinh\eta\Sin^2\chi   \cosh\eta \cos\phi
 +2\fr {g}{Bq}zy t   \sinh\eta\Sin^2\chi   \cosh\eta \sin\phi
  $$

\ses

$$
-2  \sinh^2\eta\Sin^2\chi \cos\phi    \sin\phi
\fr {g}{Bq}zxy
+
 \sinh^2\eta\Sin^2\chi \sin^2\phi
\Bigl( -  1 - \fr {g}{Bq} zy^2   \Bigr)
 $$

\ses

$$
=
\cosh^2\eta\Sin^2\chi
\Bigl(   1-\fr {g}{Bq}z t^2   \Bigr)
+2\fr {g}{Bq}z t K\cosh^2\eta\Sin^2\chi   \sinh\eta\Sin\chi
  $$

\ses

$$
+  \sinh^2\eta\Sin^2\chi
\Bigl(  - 1 - \fr {g}{Bq}zK^2\cosh^2\eta\Sin^2\chi   \Bigr),
$$
\ses
so that
$$
\fr1{J^2K^2}A_{11}=  \Sin^2\chi,
$$
or
\be
A_{11}=   K^2\fr1{h^2}\sin^2f.
\ee

{

After that, we consider the component
$A_{22}$,
getting
$$
\fr1{J^2K^2}A_{22}=
  \cosh^2\eta\Sin^2\chi \sin^2\phi
\Bigl(  - 1 - \fr {g}{Bq}zx^2   \Bigr)
$$

\ses

$$
+2  \cosh^2\eta\Sin^2\chi \cos\phi    \sin\phi
\fr {g}{Bq}zxy
+
 \cosh^2\eta\Sin^2\chi \cos^2\phi
\Bigl( -  1 -
 \fr {g}{Bq} zy^2   \Bigr).
 $$
On reducing, we are left with
\ses\\
$$
\fr1{J^2K^2}A_{22}= -\cosh^2\eta\Sin^2\chi,
$$
which results in
\be
A_{22}=  - K^2\cosh^2 \eta \fr1{h^2}\sin^2f.
\ee
Finally, we evaluate the component
$
A_{33}$,
 obtaining
 $$
\fr1{J^2K^2}A_{33}=
\sinh^2\eta\Cos^*\chi \Cos^*\chi
\Bigl(   1-\fr {g}{Bq}z t^2   \Bigr)
 $$

\ses

$$
+2\fr {g}{Bq}zx t \cosh\eta\Cos^*\chi \Cos^*\chi   \sinh\eta \cos\phi
 +2\fr {g}{Bq}zy t \cosh\eta\Cos^*\chi \Cos^*\chi   \sinh\eta \sin\phi
  $$

\ses

$$
+2 \sinh\eta\Cos^*\chi \Sin\chi
\fr {g}{Bq}q^2t
+  \cosh^2\eta\Cos^*\chi \Cos^*\chi   \cos^2\phi
\Bigl(  - 1 - \fr {g}{Bq}zx^2   \Bigr)
$$

\ses

$$
-2 \cosh^2\eta\Cos^*\chi \Cos^*\chi   \cos\phi    \sin\phi
\fr {g}{Bq}zxy
$$

\ses

$$
-2 \cosh\eta\Sin\chi \Cos^*\chi   \cos\phi
   \fr {gq}{B}    x
+
 \cosh^2\eta   \Cos^*\chi  \Cos^*\chi  \sin^2\phi
\Bigl( -  1 - \fr {gzy^2 }{Bq}   \Bigr)
 $$

\ses

$$
-2 \cosh\eta\Cos^*\chi \Sin\chi      \sin\phi
\fr {g}{Bq}  q^2  y
+
    \Sin^2\chi
       \Biggl[
-1+\fr {g} {Bq}  q^2(b+gq)\Biggr].
 $$

\ses

{

\nin
Canceling similar terms leads to
$$
\fr1{J^2K^2}A_{33}=
\sinh^2\eta\Cos^*\chi \Cos^*\chi
\Bigl(   1-\fr {g}{Bq}z t^2   \Bigr)
$$

\ses

$$
+2\fr {g}{Bq}z t K\cosh^2\eta\Cos^*\chi \Cos^*\chi   \sinh\eta \Sin\chi
+2 \sinh\eta\Cos^*\chi \Sin\chi
\fr {g}{Bq}q^2t
  $$

\ses

$$
-  \cosh^2\eta\Cos^*\chi \Cos^*\chi
-  \cosh^2\eta\Cos^*\chi \Cos^*\chi
 \fr {g}{Bq}zK^2\cosh^2 \eta\Sin\chi \Sin\chi
$$

\ses

$$
-2 \cosh\eta\Cos^*\chi \Sin\chi     \fr {g}{Bq}  q^2
   K\cosh\eta \Sin\chi
+
    \Sin^2\chi
       \Biggl[
-1+\fr {g} {Bq} q^2(b+gq)
   \Biggr]
 $$

\ses

$$
=
-\Cos^*\chi \Cos^*\chi
\Bigl(   1+\fr {g}{Bq}zK^2 \Sin\chi \Sin\chi  \Bigr)
$$
\ses

$$
-2\fr {g}{Bq}q^2K \Cos^*\chi \Sin^2\chi
+
    \Sin^2\chi
       \Biggl[
-1+\fr {g} {Bq}  q^2(b+gq)
 \Biggr].
 $$
\ses
We are left with
$$
\fr1{J^2K^2}A_{33}=
-\Cos^*\chi \Cos^*\chi
-\Cos^*\chi \Cos^*\chi
\fr{gbq}B
$$
\ses

$$
-2 \Cos^*\chi  \Sin\chi
\fr {g}{B} q^2
-      \Sin^2\chi
+  \Sin\chi         \fr {g} {B} \fr1K q^2(b+gq)
 $$

\ses

$$
=
-\Sin^2\chi    +g\Sin\chi \Cos^*\chi     -\Cos^*\chi \Cos^*\chi=-\fr 1{J^2}
$$
(we have used the identity (A.42)),
so that
\be
A_{33}=-K^2.
\ee
Thus all the coefficients in the representation (A.54) of $(ds)^2$ are valid.

Performing the choice
$z^0=\e^{\si}$ in  (A.54) yields
\be
(ds)^2= -    e^{2\si}
\Biggl[
\fr1{h^2}\sin^2 f
\Bigl((d\eta)^2-\cosh^2\eta\,(d\phi)^2\Bigr)
-(d\chi)^2
-(d\si)^2
\Biggr].
\ee
\ses
The subsequent use of the coordinates
\be
\rho=\e^{h\si}\sin(h\chi),
\qquad
\tau=\e^{h\si}\cos(h\chi)
\ee
leads to the conclusion that the metric is
 {\it conformally flat}
\be
(ds)^2= \Ka  \Bigl((ds)^2\Bigr){\Bigl|\Bigr.}_{\text{pseudo-Euclidean}},
\ee
where
\be
\Bigl((ds)^2\Bigr){\Bigl|\Bigr.}_{\text{pseudo-Euclidean}}=
-\rho^2
\Bigl((d\eta)^2-\cosh^2\eta\,(d\phi)^2\Bigr)
+(d\rho)^2+(d\tau)^2
\ee
and
\be
\Ka=\fr1{h^2}(\rho^2+\tau^2)^{(1-h)/h},
\ee
which can be written as
\be
\Ka=\fr1{h^2}\lf(K^2\rg)^{1-h}.
\ee

{

\ses

\ses

\setcounter{equation}{0}

\nin
{ \bf Appendix B:  Evaluation of spray coefficients of the
${\mathbf\cA\cR}$-space
}

\ses\ses

Let us represent the spray coefficients
$G^i=\ga^i{}_{kj}y^ky^j$
in the way
\be
G^i=g^{ik}G_k
\ee
such that
\be
G_k=y^m \D{y_k}{ x^m} - \fr12\D {F^2}{x^k}.
\ee
Using the notation
 $b_{j,k}=\partial b_j/\partial x^k$
 and
 $s_i=y^kb_{k,i}$,
we shall apply the equalities
\be
\D b{x^i}=s_i, \qquad \D q{x^i}=\epsilon \fr bqs_i+ \De,
\ee
and their implications
\be
\D {J^2}{x^k}=- \fr gB J^2 \lf(q-\epsilon \fr{b^2}q\rg)s_k
+\fr{\partial J^2}{\partial g}\fr{\partial g}{\partial x^k}
+ \De
\ee
(see (2.23)),
\be
\D{B}{x^k}=-g\fr1q(q^2+\epsilon b^2)s_k
-bq\fr{\partial g}{\partial x^k}+\De
\ee
(see (2.9)),
and
\be
\D{F^2}{x^k}=-2gq J^2s_k+
\fr{\partial F^2}{\partial g}\fr{\partial g}{\partial x^k}
+
\De
\ee
(see (2.13)),
where
 $\De$ symbolizes the summary of the terms which involve partial derivatives
of the input pseudo-Riemannian metric tensor $a_{ij}$
with respect to the coordinate variables $x^k$.

{

Taking into account the representation (2.37),  from (B.2)
we get
$$
G_k=
 -g\lf(\epsilon \fr bq(ys)b_k+qb_{km}y^m\rg)J^2
-y_k\fr gB\lf(q-\epsilon \fr{b^2}q\rg)(ys)
+
gqJ^2s_k
$$
\ses

\be
+\fr{\partial y_{k}}{\partial g} (yg)
-\fr12\fr{\partial F^2}{\partial g}\fr{\partial g}{\partial x^k}
+ \De,
\ee
where
\be
(yg)=\fr{\partial g}{\partial x^i}y^i.
\ee

To raise the index, it is convenient to apply the rules
$$
g^{ij}b_j=\fr1{F^2}\Biggl[(B+gbq)b^i
-\fr {\epsilon g}{\nu}\lf(c^2B+b(b+gc^2q)\rg)v^i\Biggr]
$$
\ses
and,
for any co-vector $t_j$,  from (2.31) we obtain
$$
g^{ij}t_j=
\Biggl[Ba^{ij}t_j
+gq(yt)b^i+\fr {\epsilon g}{\nu}
\Bigl(B(bt)-(b+gc^2q)(yt)\Bigr)v^i
\Biggr]\fr 1{F^2},
$$
where
$(yt)=y^jt_j$  and  $(bt)=b^jt_j$.
In this way we obtain
$$
G^i=-g\fr {\epsilon b}{qB}(ys)
\Biggl[(B+gbq)b^i
-\fr {\epsilon g}{\nu}\lf(c^2B+b(b+gc^2q)\rg)v^i\Biggr]
$$

\ses

$$
 -g\fr q{B}
   \Bigl[Ba^{ij}b_{jh}y^h +gq(ys)b^i
\Bigr]
 -g\fr q{B}          \fr {\epsilon g}{\nu}   \Bigl(B
 (b^jb_{jh}y^h)
 -(b+gc^2q)(ys)\Bigr)v^i
 $$

 \ses

 $$
-v^i\fr gB\lf(q-\epsilon \fr{b^2}q\rg)(ys)
+bb^i\fr gB\lf(q-\epsilon \fr{b^2}q\rg)(ys)
$$

\ses

$$
+\fr{gq}B
\Biggl[Ba^{ij}s_j
+gq(ys)b^i+\fr {\epsilon g}{\nu}
\Bigl(B(bs)-(b+gc^2q)(ys)\Bigr)v^i
\Biggr]
+E^i
+ a^i{}_{nm}y^ny^m,
$$

{

\nin
or
$$
G^i=g\fr b{qB}(ys)
\fr { g}{\nu}\lf(c^2B+b(b+gc^2q)\rg)v^i
$$

\ses

$$
 -gq
a^{ij}b_{jh}y^h
 -g\fr q{B}          \fr {\epsilon g}{\nu}   \Bigl(B
 (b^jb_{jh}y^h)
 -(b+gc^2q)(ys)\Bigr)v^i
-v^i\fr {\epsilon g}B\fr{B+gbq}q(ys)
$$

\ses

$$
+\fr{gq}B
\Biggl[Ba^{ij}s_j
+\fr {\epsilon g}{\nu}
\Bigl(B(bs)-(b+gc^2q)(ys)\Bigr)v^i
\Biggr]
+E^i
+ a^i{}_{nm}y^ny^m.
$$
Due simplifying yields
$$
G^i=g\fr b{q}(ys)
\fr { g}{\nu}(c^2-1)v^i
 -gq
a^{ij}b_{jh}y^h
 -g q          \fr {\epsilon g}{\nu}
 (b^jb_{jh}y^h)
v^i
-v^i\fr{\epsilon g}  q(ys)
$$

\ses

$$
+gq
\Biggl[a^{ij}s_j
+\fr {\epsilon g}{\nu}
(bs)v^i
\Biggr]
+E^i
+ a^i{}_{nm}y^ny^m,
$$
which can be written as
$$
G^i=  -\epsilon  \fr{g}{\nu}(ys)v^i
+
gq
\lf[a^{ij}s_j
+\fr {\epsilon g}{\nu}
(bs)v^i
\rg]
-gq
\lf[a^{ij}b_{jh}y^h
+\fr{\epsilon  g}  {\nu}
(b^jb_{jh}y^h)v^i
\rg]
$$

\ses

$$
       +E^i
+a^i{}_{nm}y^ny^m.
$$
This method results in the representations (2.61)-(2.63).

{

\ses

\ses

\setcounter{equation}{0}

\nin
{ \bf Appendix C:  Conformal transformation in the
$ {\mathbf\cU\cA\cR}$-space
}

\ses

\ses

It proves possible to indicate
the explicit transformation
\be
\zeta^m=\zeta^m (x,y)
\ee
to fulfill the conformal claim.
Indeed, let us take
\be
\zeta^m=\Bigl[hv^m  -   ( b+\fr12gq)b^m    \Bigr]J\fr1{\varkappa h},
\qquad   \varkappa=\fr1{h}\lf(|F^2|\rg)^{(1-h)/2},
\ee
where $v^m=y^m+bb^m$ and both the values
$\epsilon=1$ or $\epsilon=-1$ are admissible.
Evaluating the derivatives
\be
\zeta^m_n~:=\D{\zeta^m}{y^n}
\ee
yields
\be
\zeta^m_n
=E^m_n
+\fr1N \zeta^mC_n
-
\fr1{\ka}\ka_n\zeta^m,
\ee
where
\be
E^m_n
= \Biggl[   h(\de^m_n+b_nb^m) -  \lf(  b_n+\epsilon   \fr1{2q}gv_n\rg)b^m\Biggr]
J\fr1{\varkappa h}
\ee
\ses
and
\be
\ka_n=\D{\ka}{y^n}=(1-h)y_n\fr1{F^2}\ka,
\ee
and we have used the equality
\be
\D J{y^n} =  \fr1NC_n
\ee
(which is a direct implication of the formulas (2.32) and (2.41)).
It is useful to  note that
\be
E^m_ny^n=\zeta^m , \qquad    E^m_nb^n=   b^m J\fr1{\varkappa h}.
\ee

{

Let us find the transform
\be
s^{ij}~:= g^{mn}\zeta_m^i\zeta_n^j
\ee
of the initial Finslerian metric tensor
$g^{mn}$. To this end we apply
the representation (2.39) of $g^{mn}$ (with $c=1$ which entails $\nu=q$)
and note the vanishing
$$
-
\fr{\zeta^i\zeta^j}{N^2}
C_mC^m
-\fr{(1-h)^2}{F^2}\zeta^i\zeta^j
+\fr {1-h}{F^2}
(\zeta^iy^m E_m^j+\zeta^jy^m E_m^i)=0,
$$
getting
$$
s^{ij}
=
\fr1{N}(\zeta^iC^m E_m^j+\zeta^jC^m E_m^i)
$$

\ses

$$
+
\Biggl[
a^{mn}E_m^iE_n^j
+ \epsilon\fr gqbb^mb^nE_m^iE_n^j
+ \epsilon\fr gq  (b^my^n+y^mb^n)    E_m^iE_n^j
- \epsilon\fr g{Bq}(b+gq)  y^m  y^n E_m^i  E_n^j
\Biggr]
\fr1{J^2}.
$$
Simple calculation shows that
$$
 a^{mn}E_m^iE_n^j=
h a^{in}
\Biggl[
   h(\de^j_n+b_nb^j) +   \lf(    b_n+\epsilon \fr1{2q}gv_n\rg)b^j
\Biggr] \fr{J^2}{\ka^2h^2}
+
b^i
\Bigl[
(h+1)b^n+\epsilon \fr1{2q}gv^n\Bigr]\fr{E_n^jJ} {\varkappa h}
$$

\ses

$$
=
h
\Biggl[
   h(a^{ij}+b^ib^j)  +    \lf(    b^i+\epsilon \fr1{2q}gv^i\rg)b^j
\Biggr]\fr{J^2} {\ka^2h^2}
-
\Biggl[
(h+1)+\epsilon \fr1{2q}gb\Biggr] \fr{b^ib^jJ^2} {\ka^2h^2}
+\epsilon  \fr1{2q}g\zeta^j \fr{Jb^i}{\varkappa h},
$$
\ses
so that
$$
 a^{mn}E_m^iE_n^j=
h
\Biggl[
   h(a^{ij}+b^ib^j)   + \epsilon   \fr1{2q}gv^i  b^j
\Biggr] \fr{J^2}{\ka^2h^2}
-
\Biggl[
1+\epsilon \fr1{2q}gb\Biggr]b^ib^j
\fr{J^2}{\ka^2h^2}
+b^i\epsilon  \fr1{2q}g\zeta^j \fr J{\varkappa h}.
$$

 {

In this way we get
$$
h^2s^{ij}
=
\epsilon
\fr g{qF^2}(b+gq) \zeta^i\zeta^j
h^2
+\epsilon
\fr g2\fr1{q}
(\zeta^i b^j+\zeta^jb^i)  J\fr1{\varkappa} h
$$

\ses

$$
+h^2a^{ij}\fr1{\ka^2}
   +h^2b^ib^j\fr1{\ka^2}
+\epsilon
\fr1{2q}gb^j hv^i\fr1{\ka^2}
-
\Biggl[
1+\epsilon
\fr1{2q}gb\Biggr]b^ib^j
\fr1{\ka^2}
+b^i\epsilon  \fr1{2q}g\zeta^j   J\fr1{\varkappa} h
$$

\ses

$$
+\epsilon
 \fr gqb b^ib^j\fr1{\ka^2}
-\epsilon
 \fr gq  (b^i\zeta^j+\zeta^ib^j)    J\fr1{\varkappa} h
-\epsilon
 \fr g{Bq}(b+gq) \zeta^i\zeta^j
\fr1{J^2}
h^2,
$$
\ses
or
$$
h^2s^{ij}
=
\epsilon
\fr g2\fr1{q}
(\zeta^i b^j+\zeta^jb^i)    J\fr1{\varkappa} h
+h^2a^{ij}\fr1{\ka^2}
   +h^2b^ib^j\fr1{\ka^2}
+\epsilon
\fr1{2q}gb^j hv^i\fr1{\ka^2}
$$

\ses

$$
-
\Biggl[
1+\epsilon
\fr1{2q}gb\Biggr]b^ib^j
\fr1{\ka^2}
+b^i\epsilon  \fr1{2q}g\zeta^j    J\fr1{\varkappa} h
+\epsilon
 \fr gqb b^ib^j\fr1{\ka^2}
-\epsilon
 \fr gq  (b^i\zeta^j+\zeta^ib^j)  J\fr1{\varkappa} h.
$$
\ses
Inserting  here the equality
$$
hv^i=\ka
h\fr1J\zeta^i-
\lf(b+\fr12gq\rg)b^i
$$
(which ensues from (C.2)),
we eventually  find
\be
s^{ij}=   \fr1{\ka^2} a^{ij}.
\ee

 From (C.4), (C.10), and (2.32)    it follows that
\be
\det(\zeta^i_m)   =(J /  \varkappa)^N >0 .
\ee
The metric tensor transformation (C.9) can be inverted to read
\be
g_{mn}=\varkappa^2 a_{ij}\zeta^i_m \zeta^j_n.
\ee

The functions (C.2)  possess the property  of homogeneity of degree
 $h$:
\be
\zeta^i(x,ky) =
k^h\zeta^i(x,y), \qquad k>0,\,\forall y,
\ee
which entails the identity
\be
y^n\zeta^i_n \equiv h\zeta^i.
\ee
Therefore, from (C.12) and the equality $g_{mn}y^my^n=F^2$ we can infer that
\be
F^2=h^2\varkappa^2    S^2(x,\zeta),
\ee
\ses
where
\be
S^2(x,\zeta)=a_{ij}(x)\zeta^i\zeta^j.
\ee
The equality (C.12) motivates introducing the tensor
\be
t_{mn}(x,\zeta)=\varkappa^2 a_{mn}.
\ee

From (C.2) and (C.15) we can obtain the remarkable equality
\be
|F^2(x,y)|^{h(x)}=|S^2(x,\zeta)|.
\ee
So, if we introduce the scalar
\be
p(x,\zeta) =\varkappa^2,
\ee
from (C.2)  we  get the representation
\be
p(x,\zeta)=\fr1{h^2(x)}|S^2(x,\zeta)|^{(1-h(x))/h(x)}
\ee
which enables us to read
 the tensor (C.17) as follows:
\be
t_{mn}(x,\zeta)=p(x,\zeta) a_{mn}(x).
\ee

This line of reasoning leads to  the definition
\be
{\mathbf\cC}_{g} :=
\{\cR_{N};\,TM;\, \zeta\in TM;\,\,g(x);\,t_{mn}(x,\zeta)\}
\ee
which we  call  the {\it factor-pseudo-Riemannian space.}
The transformation (C.1)-(C.4)  assigns the map
\be
{\cal Z}:~ {\mathbf\cU\cA\cR}_g
\to {\mathbf\cC_g},
\ee
which sends vectors to vectors
\be
y\stackrel{\cal Z}{\to} \zeta
\ee
(at each point $x\in M$)
and replaces  the Finslerian metric tensor
$\{g_{mn}(x,y)\}$
by  the factored tensor (C.21):
\be
\lf\{g_{mn}\rg\}\stackrel{\cal Z}{\to} \lf\{t_{mn}\rg\}.
\ee
Obviously, the map is the diffeomorphism.

Thus we are entitled to formulate the following proposition.

\ses

\ses

PROPOSITION. {\it  The }  $ {\mathbf\cU\cA\cR}$-space }
{\it is
metrically isomorphic to the factor-pseudo-Riemannian space}:
\be
g_{mn}(x,y)dy^mdy^n=t_{mn}(x,\zeta)d\zeta^md\zeta^n.
\ee

\ses

\ses

Since $v^i=q=0$ at $y^i=b^i$, from (C.2) it follows that
\be
\zeta^i(x,b)=b^i(x)
\ee
(the equality $K^2(x,b)=1$ has been into account;
examine  (4.10)-(4.14)).
Therefore, we can formulate the following.

\ses

\ses

PROPOSITION. {\it  The involved preferred vector field $b^i(x)$ is the proper element
of the  $ {\mathbf\cU\cA\cR}$-space  isomorphism}:
\be
\lf\{b^i(x)\rg\}\stackrel{\cal Z}{\to} \lf\{b^i(x)\rg\}.
\ee

\ses

\ses

The transformation (C.1) is trivial in the pseudo-Riemannian limit:
$$
\zeta^i\bigl|_{g=0}\bigr.=y^i.
$$

Let us examine the indicatrix curvature.
At each point $x$,
the invention  (C.23)-(C.25) maps isomorphically
 the Finslerian indicatrix (defined by $|F^2(x,y)|=1$)
onto the pseudo-Euclidean sphere $ \cS_x({\mathbf\cC_g})$
(defined by  $|S^2(x,\zeta)|=1$) in the space $ {\mathbf\cC_g}$.
The factor $p$, having been proposed by (C.20), equals $1/h^2(x)$
 on this $ \cS_x({\mathbf\cC_g})$  and thereby
 influences the form of the metric tensor
induced on $ \cS_x({\mathbf\cC_g})$.
Namely, let some coordinates $m^a$  be introduced on
$ \cS_x({\mathbf\cC_g})$, so that the tensor is of the type
$ i_{ab}= i_{ab}(x,m^c)$.
For instance, in the time-like sector we can take
$   m^a=\zeta^a/\zeta^0, \, \epsilon =1.     $
Let $L^i$ denote  the unit vectors,  so  that $\epsilon L^iL^jt_{ij}=1$.
We use the coordinates to parameterize the vectors, obtain the projection factors
$L^i_a=\partial L^i/\partial m^a$,
and construct the tensor
$
 i_{ab}=t_{ij}L^i_a L^j_b.
$
In view of the factored structure (C.20)-(C.21) of the tensor $t_{ij}$,
we get  the equality
$$
  i_{ab}=\fr1{h^2}
 \widehat i_{ab},
$$
 where
$
  \widehat i_{ab}=a_{ij}L^i_a L^j_b
$
   is the indicatrix metric tensor
that  is obtainable  in the pseudo-Riemannian geometry proper.
Since $h$ is independent of $\zeta$,
 the associated Christoffel symbols
$$
i_a{}^c{}_b
=\fr12i^{ce}\lf(\D{i_{ea}}{m^b}+\D{i_{eb}}{m^a}-\D{i_{ab}}{m^e}\rg),
\qquad
 \widehat i_a{}^c{}_b
=\fr12 \widehat   i^{ce}\lf(\D{ \widehat i_{ea}}{m^b}+\D{ \widehat i_{eb}}{m^a}
-\D{ \widehat i_{ab}}{m^e}\rg)
$$
are equivalent:
$$
  i_a{}^c{}_b   =   \widehat i_a{}^c{}_b .
$$
Therefore,  the indicatrix curvature tensor
$$
I_a{}^c{}_{bd}=
\D{i_a{}^c{}_b}{m^d}-
\D{i_a{}^c{}_d}{m^b}
+i_a{}^e{}_b    i_e{}^c{}_d
-i_a{}^e{}_d    i_e{}^c{}_b
$$
 is identical to  the tensor
$\widehat I_a{}^c{}_{bd}$
constructible by the same rule from the tensor
$ \widehat i_{ab}$,
that is,
$I_a{}^c{}_{bd}
= \widehat I_a{}^c{}_{bd}$.
Let us now consider the tensors
$
I_{acbd}  = i_{ce}I_a{}^e{}_{bd}
$ and
$\widehat I_{acbd}
= \widehat i_{ce}\widehat I_a{}^e{}_{bd}.
$
Since
$
\widehat I_{acbd}=  -
\epsilon
(\widehat i_{cb}\widehat i_{ad}-\widehat i_{cd}\widehat i_{ab})
$
is the ordinary case characteristic of  the pseudo-Riemannian geometry
(which  reflects the fact that  the indicatrix  constant curvature is of the unit type),
 we get
$
I_{acbd}=  -
\epsilon (\widehat i_{cb}\widehat i_{ad}-\widehat i_{cd}\widehat i_{ab})/h^2,
$
which in turn entails
\be
I_{acbd}=  -
\epsilon  h^2
(i_{cb}i_{ad}-i_{cd}i_{ab}),
\ee
whence
\be
 \cR_{\text{${\mathbf\cU\cA\cR}$-pseudo-Finsleroid  indicatrix} }=- \epsilon h^2.
\ee

 Amazingly,    the angle (2.71)-(2.72)   may  be written as
  \be
 \al_{\{x\}}(y_1,y_2) =  \fr1h
 \al_{\text{pseudo-Riemannian}},
 \ee
 where
 \be
  \al_{\text{pseudo-Riemannian}}= \arccosh
  \fr{a_{mn}(x)\zeta^m (x,y_1) \zeta^n(x,y_2)}
 { \sqrt{S^2(\zeta (x,y_1))}\,\sqrt{S^2(\zeta (x,y_2))} }, \,\, \epsilon =1,
\ee
\ses
and
\be
 \al_{\text{pseudo-Riemannian}}= \arccos
  \fr{-a_{mn}(x)\zeta^m (x,y_1) \zeta^n(x,y_2)}
 { \sqrt{\lf|S^2(\zeta (x,y_1))\rg|}\,\sqrt{\lf|S^2(\zeta (x,y_2))\rg|} }, \,\, \epsilon =-1,
  \ee
(insert here (C.2) and make comparison with  (2.71)-(2.73)).
The conformal factor $p$ does not enter the right-hand parts of (C.32)-(C.33).
The angle  (C.32)-(C.33) belongs to the sense of the ordinary
pseudo-Riemannian geometry and operates in the
factor-pseudo-Riemannian space (C.22).

{

If we insert (C.2) in (C.16), we obtain the useful
equality
\be
\lf( \fr J{\varkappa h}\rg)^2=\fr{S^2(x,\zeta)}{B(x,y)}.
\ee

Also, from (C.2) it follows that
\be
\zeta^ib_i=\lf(b+\fr12gq \rg) \fr J{\varkappa h}
\ee
and
\be
\zeta^i+ (\zeta^nb_n)b^i=
 \fr J{\varkappa h}
 hv^i.
 \ee
According to
(2.86) and (2.91), we have
\be
\chi   = \fr1h    \arcsinh \fr{\zeta^nb_n(x)}    {\sqrt{ S^2(x,\zeta) }},\, \, \epsilon =1;
\qquad
\chi   = \fr1h    \arccos \fr{\zeta^nb_n(x)}    {\sqrt{| S^2(x,\zeta)|}}, \,\, \epsilon =-1,
 \ee
and from (2.84) and (2.89) we get
\be
b=F\Sh\chi, \quad \epsilon =1; \qquad  b=K\Cos\chi, \quad \epsilon =-1.
\ee
From (C.19)-(C.20) and (2.78)
it follows that
\be
h\varkappa=
|S^2(x,\zeta)|^{(1-h)/(2h)}, \qquad
\fr1J=e^{\frac12g\chi}.
\ee

The indicated formulas allow us to write down the explicit form of the
inverse
to (C.1) and (C.2), namely we find
\be
y^i=y^i(x,\zeta)
\ee
with
\be
y^i=-bb^i+
\fr1h
\Bigl(\zeta^i+ (\zeta^nb_n)b^i\Bigr)
 \fr {h\varkappa }J
 \ee
(examine (C.36) and remind that $v^i=y^i+bb^i$).

It is possible to find straightforwardly the coefficients
 \be
y^i_j~:= \D{y^i}{\zeta^j} .
 \ee

{

\ses

\ses

\def\bibit[#1]#2\par{\rm\noindent\parskip1pt
                     \parbox[t]{.05\textwidth}{\mbox{}\hfill[#1]}\hfill
                     \parbox[t]{.925\textwidth}{\baselineskip11pt#2}\par}

\nin {  REFERENCES}

\ses

\bibit[1] H. Rund: \it The Differential Geometry of Finsler
 Spaces, \rm Springer, Berlin 1959.

\bibit[2] D. Bao, S. S. Chern, and Z. Shen: {\it  An
Introduction to Riemann-Finsler Geometry,}  Springer, N.Y., Berlin
2000.

\bibit[3] J. I. Horv\'ath: New geometrical methods of the theory of physical fields,
\it Nouv. Cim. \bf 9 \rm(1958), 444--496.

\bibit[4] Y. Takano: Theory of fields in Finsler spaces,
 \it Progr. Theor. Phys. \bf32 \rm(1968),  1159--1180.

\bibit[5] R. Mrugala: Riemannian and Finslerian geometry in thermodynamics,
\it Open Syst. Inform. Dyn. \bf1 \rm    (1992), 379--396.

\bibit[6] R. S. Ingarden and L. Tamassy:
 On Parabolic geometry and irreversible macroscopic time,
 \it Rep. Math. Phys. \bf 32 \rm(1993), 11.

\bibit[7] R. S. Ingarden:  On physical applications of Finsler  Geometry,
   \it Contemporary Mathematics \bf 196 \rm(1996), 213--223.

\bibit[8]    G. S. Asanov:   Finsleroid--Finsler  space with Berwald and  Landsberg conditions,
  arXiv:math.DG/0603472 (2006);
 Finsleroid--Finsler spaces of positive--definite and relativistic types,
\it Rep. Math. Phys.  \bf 58 \rm(2006), 275-300.

\bibit[9]    G. S. Asanov:    Finsleroid-Finsler  space and spray   coefficients,
  arXiv:math.DG/0604526 (2006);
  Finsleroid--Finsler space and geodesic spray    coefficients,
{\it Publ.  Math. Debrecen } {\bf 71/3-4} (2007), 397-412.

\bibit[10] G. S. Asanov:   Finsleroid-regular    space  developed.    Berwald case,
  arXiv:math.DG/0711.4180  (2007);
 Finsleroid-regular    space:~ curvature tensor,
continuation of gravitational Schwarzschild metric,
arXiv:math-ph/0712.0440  (2007);
    Finsleroid-regular    space.  Gravitational  metric.
  Berwald case, \it Rep. Math. Phys.
 \bf 62 \rm(2008), 103-128;
 Finsleroid-regular   space.    Landsberg-to-Berwald implication,
  arXiv:math.DG/0801.4608  (2008);
   Finsleroid  corrects     pressure  and energy   of
  universe. Respective  cosmological equations,
 arXiv:math-ph/0707.3305 (2007);
Finsleroid-cosmological equations,   \it Rep. Math. Phys. \bf 61 \rm(2008), 39-63.

\bibit[11]   G. S. Asanov:
Finslerian anisotropic relativistic metric function obtainable under breakdown of
rotational symmetry,
  arXiv:gr-qc/ 0204070  (2002);
Finslerian post-Lorentzian kinematic transformations in anisotropic-space time,
  arXiv:gr-qc/0207117  (2002).

\bibit[12]   G. S. Asanov:   Finsleroid--Finsler  space of involutive case,
  arXiv:math.DG/0710.3814  (2007);
 Finsleroid--Finsler  space of involutive case and $A$-special relation,
  {\it Publ.  Math. Debrecen } {\bf 72/3-4} (2008), 737-747.

\end{document}